\magnification=\magstep1
\input amstex
\documentstyle{amsppt}
\catcode`\@=11 \loadmathfont{rsfs}
\def\mycal{\mathfont@\rsfs}
\csname rsfs \endcsname \catcode`\@=\active  

\vsize=7.5in

\topmatter 
\title On ergodic embeddings of factors  \\ 
$\text{\it To Dick Kadison, in memoriam}$ \endtitle

\author  Sorin Popa \endauthor
\affil     {\it  University of California, Los Angeles} \endaffil

\rightheadtext{Ergodic embeddings of factors}

\address Math.Dept., UCLA, Los Angeles, CA 90095-1555, USA \endaddress
\email  popa\@math.ucla.edu\endemail

\thanks Supported in part by NSF Grant DMS-1700344 and the Takesaki Chair in Operator Algebras \endthanks

\abstract An inclusion of von Neumann factors $M \subset \Cal M$  is {\it ergodic} 
if it satisfies the irreducibility condition $M'\cap \Cal M=\Bbb C$. We investigate  
the relation between this and several stronger ergodicity properties, such as $R$-{\it ergodicity}, 
which requires $M$ to admit an embedding of the hyperfinite II$_1$ factor $R\hookrightarrow M$ 
that's ergodic in $\Cal M$. We prove that if $M$ is {\it continuous} (i.e., non type I)  
and contains a maximal abelian $^*$-subalgebra 
of $\Cal M$, then $M\subset \Cal M$ is $R$-ergodic. This shows in particular that 
any continuous factor contains an ergodic copy of $R$.  
\endabstract 

\endtopmatter

\document

\heading 0. Some background motivation \endheading

Dick Kadison liked to recount of how John von Neumann, in search of a good problem to give to Fred Murray,  
his freshly assigned postdoc at Princeton  
in the mid 1930s, something that would both interest himself 
and fit Murray's background, opted to follow up on his 
earlier work [vN29], on ``rings of operators''. 
The concrete project he proposed was  to investigate whether there exist any other factors acting on a separable Hilbert space besides the {\it atomic} ones 
$\Cal B(\ell^2_n)$, $1\leq n \leq \infty$ (later called of {\it type} $I_n$).  Kadison 
enjoyed watching his younger interlocutor  being taken aback, for a brief moment, at the reckoning that there was indeed a time 
when it was not clear at all what the answer to this question could be !

The ensuing discovery of  {\it continuous} factors in [MvN36], namely factors that have no minimal projections (atoms), 
which Murray-von Neumann  labelled of type II and III, according to the possible dimension function they  admit, 
unravelled a series of striking new phenomena. It is well known that von Neumann was fascinated by the 
existence of  the {\it continuous dimension} on type II factors (ranging over $[0,1]$ in the II$_1$ case, and over $[0, \infty]$ in the II$_\infty$ case), 
and the new perspectives this could bring to various fields of mathematics (cf. also [vN54]).  

But he and Murray seemed equally enthralled by the fact that, unlike the type I case where 
any von Neumann subalgebra $M \subset \Cal M=\Cal B(\ell^2)$ satisfies the bicommutant property $(M' \cap \Cal M)'\cap \Cal M=M$ ([vN29]), 
for certain continuous factors $\Cal M$ of type II and III  this may fail to be true 
(see pages 185, 209, 229 in [MvN36]). The examples found 
in ([MvN36]) are inclusions of  factors 
obtained via the {\it group measure space} (crossed-product type) construction:  let $\Gamma \curvearrowright (X, \mu)$ 
be a free, ergodic action by non-singular 
transformations of a countable group $\Gamma$ on a measure space $(X, \mu)$, take $\Cal M=L^\infty(X) \rtimes \Gamma$ as ambient factor and 
$M=L^\infty(X)\rtimes H \subset \Cal M$ as subfactor, for some $H\subset \Gamma$; if $\Gamma \neq H\neq 1$, then the bicommutant property fails, 
in fact in a rather ``dramatic way'', 
as one even has $M'\cap \Cal M=\Bbb C$ whenever $H$ itself acts ergodically. Due to their subsequent work ([MvN43]), 
if one takes $\Gamma=S_\infty$, $H$ a proper infinite subgroup and the action to be Bernoulli, this example covers the case $\Cal M$ is the 
(unique) {\it approximately finite dimensional} (AFD) II$_1$ factor $R:=\overline{\otimes}_n (\Bbb M_{2}(\Bbb C), tr)_n$,  
later called the {\it hyperfinite} II$_1$ factor in [D57], an algebra of fundamental importance. So $R$ contains subfactors $M \subsetneq R$ satisfying $M'\cap R=\Bbb C$. 

One of the problems posed in ([MvN36]) asked whether {\it any} continuous factor  $\Cal M$ may in fact 
contain subfactors $M\subset \Cal M$ failing the bicommutant property (see Problem 10 on page 185 therein; note that this was reiterated 
as Problem 10 in [K67]). For type II factors, 
this was answered in  the affirmative by Fuglede and Kadison ([FK51]), who noticed that the commutant $M'\cap \Cal M$ 
of any maximal hyperfinite subfactor $M=R$ of a II$_1$ factor $\Cal M$  is either $\Bbb C$, or non-factorial, both cases implying that if $\Cal M\not\simeq R$ then 
$(M'\cap \Cal M)'\cap M \neq M$. The type III case was settled 
more than two decades later,    
as a ``side effect'' of  the Tomita-Takesaki theory and Connes' ground-breaking work on  decomposition of 
type III factors ([C73]; cf. also [CT76]). So indeed, as predicted in ([MvN36]), the bicommutant property only holds in factors of type I. 

It  became more and more apparent over the years 
that in fact the most natural, generic position for an inclusion of continuous factors $M \subset \Cal M$ is the one satisfying the {\it irreducibility} condition  
$M'\cap \Cal M= \Bbb C1$. This amounts to $\Cal U(M)$ acting ergodically on $\Cal M$, via Ad, 
hence the terminology of {\it ergodic embedding} that we will adopt here, emphasizing the ``dynamical'' interaction between algebras. 
For instance, crossed products and amalgamated free product constructions/decompositions, 
which now play a key role in studying the structure and classification of type II and III factors, 
generally give rise to ergodic embeddings $M \subset \Cal M$. 

At the same time, the evolution of the subject generated a variety of questions about 
embedding factors one into another, where adding the ergodicity condition is of key interest. 
A special case in point is when one of the factors involved is the hyperfinite II$_1$ factor $R$.   

It was already noticed in ([MvN43]) that $R$ embeds into 
any other II$_1$ factor  $\Cal M$ (in fact, in any  infinite dimensional factor). 
Starting with the mid-60s, the question of whether $R$ is the ``smallest'' continuous factor, in that any 
II$_1$ factor embeddable into it is isomorphic to $R$, became of fundamental importance (NB: this is already 
addressed on the last lines of page 717 in [MvN43]). 
This was in close relation to  the development of the notion of {\it amenability}  for factors, 
in several equivalent ways, see [C75] for an account and the names 
involved in this extraordinary body of work. It all culminated with Connes'  famous theorem [C76]: all 
amenable II$_1$ factors (in particular all II$_1$ subfactors of $R$) are isomorphic to $R$. 

The question of whether $R$ can be embedded ergodically into any II$_1$ factor was answered in [P81a]. 
In its full generality, the result shows that if $M\subset \Cal M$ is an irreducible inclusion of separable II$_1$ factors, 
then $M$ contains an ``$R$-direction'' that's ergodic in $\Cal M$:  there exists $R\hookrightarrow M$ 
such that $R'\cap \Cal M=\Bbb C$. We will call $R$-{\it ergodicity} this strengthened form of ergodicity for an inclusion of  factors $M\subset \Cal M$. 

Moreover, [P81a] shows that if $M\subset \Cal M$ is an ergodic inclusion of II$_1$ factors, 
then one can construct the embedding $R\hookrightarrow M$ such that the natural ``diagonal'' of $R$, $D:=\overline{\otimes}_n 
(D_2, tr)_n$, is maximal abelian (a {\it MASA}) in $\Cal M$. This also answered in the affirmative the 
II$_1$ factor case of Problem 11 in [K67], about 
whether the ergodicity of an inclusion of factors $M\subset \Cal M$ entails the existence 
in $M$ of a maximal abelian $^*$-subalgebra (abbreviated {\it MASA}) of $\Cal M$. The latter property, due to Kadison, 
is easily seen to imply ergodicity for arbitrary inclusions of factors.  To be consistent with the ``dynamical'' point of view  
adopted here,  we will call it {\it MASA-ergodicity}. Let us also point out that, as shown 
by several examples in [P81b], a maximal 
hyperfinite II$_1$ subfactor of a II$_1$ factor $\Cal M$ may actually fail to be ergodic in $\Cal M$. 

Producing ``large'' ergodic copies of $R$ inside arbitrary factors, and more generally inside  
irreducible inclusions of factors $M \subset \Cal M$, in the spirit of [P81a], i.e., 
establishing $R$-ergodicity from mere ergodicity, turns out to be of crucial importance for a multitude of problems,  
notably in proving vanishing cohomology results (cf. e.g.,  [P17], [PV14]). 
This is because once having $R\subset M$ that's ergodic in some appropriate ``augmentation'' $\Cal M$ of $M$, the amenability of $R$  can be used to 
``push'' any $x\in \Cal M$ into  $R'\cap \Cal M=\Bbb C1,$ by averaging over 
unitaries in $R$, via the Ad-action (cf. [Sc63]). When applied to suitable $x$, 
this amounts to ``untwisting'' a cocycle (see e.g., [P18a], [PV14]). 

In case $M\subset \Cal M$ are II$_1$ factors and the $x$'s that need averaging are from $\Cal M$, then the Hilbert-space 
structure given by the trace on $\Cal M$ allows averaging over the entire $\Cal U(M)$ and the Hilbert-space convexity 
shows that $1$ together with one of the unitaries (thus an ``abelian direction'') already pushes $x$ ``a little bit'' into $M'\cap \Cal M$. 
It is this fact that's being used in the proof of [P81a]. 

So another strengthening of ergodicity naturally comes into picture:  the condition that any $x\in \Cal M$ can be pushed into scalars by averaging over the unitaries in $M$, 
i.e., $\overline{\text{\rm co}}^w \{uxu^* \mid u\in \Cal U(M)\} \cap \Bbb C1 \neq \emptyset$, $\forall x\in \Cal M$. 
Since it involves convex combinations of conjugates of $x$, in the spirit of 
Dixmier's averaging theorem [D57], but taking weak closure instead of norm closure and the unitaries in $M$ instead of $\Cal M$, 
this condition has been called {\it weak 
relative Dixmier property} in ([P98]). But instead, in order to emphasize the dynamic aspect of this property 
that's reminiscent of von Neumann's mean ergodic theorem, we call it here {\it MV-ergodicity} (as an abbreviation of   
{\it mean value} ergodicity). One can easily see that this condition is implied by both $R$-ergodicity and MASA-ergodicity. 

Deciding whether a specific ergodic inclusion $M\subset \Cal M$ is $R$-ergodic, MASA-ergodic, or MV-ergodic, or 
that it does not satisfy one of these conditions, is usually a difficult and subtle problem.  
Some of the most important  open questions in operator algebras can be reduced to solving this type of problem for certain $M\subset \Cal M$. 

For instance, as shown in [H87], Connes bicentralizer  conjecture for a III$_1$ factor $\Cal M$ would hold true if one could prove 
that the ergodic $\text{\rm II}_\infty \subset \text{\rm III}_1$  inclusion $M\subset \Cal M$ associated with the continuous decomposition of $\Cal M$ is MV-ergodic  
(so showing $R$-ergodicity or MASA-ergodicity would be sufficient as well).  

Also, as explained in (Section 5 of [P18]), one possible approach towards proving that the free group factors are non-isomorphic 
and that $L(\Bbb F_\infty)$ is $\infty$-generated, is to show that if a II$_1$ factor $M$ is {\it stably single generated} 
then $M$ necessarily contains an  irreducible hyperfinite subfactor $L \subset M$ such that $M \subset \Cal M=\langle M, e_L \rangle$ is 
$R$-ergodic. As we will explain in a forthcoming paper, establishing $R$-ergodicity for certain ergodic II$_1 \subset \text{\rm II}_\infty$ inclusions  is also relevant to Connes Embedding problem 
(see also Section 6 in [P18a]).  

\heading 1. Results in this paper
\endheading 

Our purpose here is to investigate the relation between the various notions of ergodicity for  embeddings of factors $M \subset \Cal M$  emphasized above: 
$R$-ergodicity, MASA-ergodicity and MV-ergodicity, accounting for the existence of an ergodic $R$-direction, ergodic abelian direction, 
and mean value type convergence, of the action $\Cal U(M) \curvearrowright^{\text{\rm Ad}} \Cal M$. 

By [P81a], for inclusions of II$_1$  factors these conditions are all equivalent to plain ergodicity. 
But as shown in [GP96], there are examples of $\text{\rm II}_1 \subset \text{\rm II}_\infty$ ergodic inclusions 
that are not $R$-ergodic, nor MASA-ergodic: if $L\subset M=L(\Bbb F_2)$ is any irreducible hyperfinite subfactor of the free group factor 
(which always exists by [P81a]), then the associated 
basic construction subfactor $M \subset \langle M, e_{L} \rangle = \Cal M$ is irreducible, with $\Cal M$ of type II$_\infty$, 
yet there exists no hyperfinite subfactor $R\subset M=L(\Bbb F_2)$ that's irreducible in $\Cal M$, and no abelian 
$A\subset M$ that's a MASA in $\Cal M$. 

While these were the only known examples of ergodic embeddings that are not MASA-ergodic nor $R$-ergodic, it was not known whether 
they fail MV-ergodicity as well. In fact, no examples of ergodic but not MV-ergodic embeddings were known. 
 
Our first result answers both of these problems, by providing examples of inclusions of factors that are ergodic but not MV-ergodic, 
and showing existence of inclusions of factors that are MV-ergodic but not $R$-ergodic. 

\proclaim{1.1. Theorem}$1^\circ$ Let $\Gamma$ be a nonamenable group and $\Gamma \curvearrowright R=\overline{\otimes}_{g\in \Gamma} (\Bbb M_2(\Bbb C), tr)_g$ 
be the Bernoulli $\Gamma$-action with base $(\Bbb M_2(\Bbb C), tr)$. Let $M=R \rtimes \Gamma$  and $R \subset M \subset \langle M, e_R \rangle \simeq R^\infty:= R\overline{\otimes} \Cal B(\ell^2\Bbb N)$,  
be the associated crossed product $\text{\rm II}_1$ factor and its basic construction inclusion. Then $M\subset R^\infty$ is ergodic but not MV-ergodic. 

$2^\circ$  Any interpolated free group factor $L(\Bbb F_t)$, $1< t \leq \infty$, admits an embedding into a properly infinite 
AFD factor  that's MV-ergodic but not $R$-ergodic. 
\endproclaim 

To prove the first part, we show that if $M\subset R^\infty=\langle M, e_R \rangle$ is MV-ergodic, then the group $\Gamma$ admits an invariant mean, 
a fact that contradicts its non-amenability. 

The second part follows by combining a recent result in [DPe19], which shows that 
the boundary inclusion of a II$_1$ factor, $M\subset \Cal B_\varphi$, associated with a Markov averaging operator $\varphi=\varphi_{\Cal U}$ corresponding to a countable 
set $\Cal U=\Cal U^* \subset \Cal U(M)$ that generates $M$, is MV-ergodic, with the result in [GP96], showing that if $M=L(\Bbb F_t)$ and $L\subset M$ is an ergodic 
hyperfinite subfactor, then $M\subset \langle M, e_L \rangle$ cannot be  $R$-ergodic. 

In fact, any example of an MV-ergodic inclusion that's not $R$-ergodic, is not MASA-ergodic either. This  is due to our second and main result in this paper, 
which shows that for any embedding of a continuous factor into another factor, 
MASA-ergodicity implies $R$-ergodicity.   

\proclaim{1.2. Theorem} Let $M\subset \Cal M$ be an embedding of continuous, separable factors. 
If $M \subset \Cal M$ is MASA-ergodic, then it is $R$-ergodic. 
\endproclaim

A key ingredient in the proof of this result is a property of MASAs in factors, $A\subset M$,  that we call {\it flatness}. It requires that 
for any finite set of unit vectors $F$ in the Hilbert space on which $M$ acts, there exists a unitary $u\in M$ such that 
the vector states implemented by $F$ ``almost coincide'' when restricted to $uAu^*$. 
In other words, one can rotate $A$ in $M$ so that all measures on $A$ implemented by 
$F$ look alike. Flatness has been  established in [P17], [P18b] for MASAs in II$_1$ factors. We prove here that in fact any MASA,  
in any continuous factor, is flat (see Theorem 5.5). 
We in fact also clarify  the atomic, type I case, where the complete opposite is true:  
no MASA of $M=\Cal B(\Cal H)$ can be flat.  

Once $M$ contains a MASA of $\Cal M$ that's flat in $M$, one can approximate it by finite dyadic partitions that are equivalent in $M$, 
have a given  $F\subset \Cal H$ implement ``almost the same''  state on it, 
and get ``closer and closer'' to being maximal abelian in $\Cal M$ (using a local characterization of MASAs inspired by [P84], [Sk77], see Lemma 4.3). We build dyadic 
matrix units $\{e_{ij}\}_{i,j} \subset M$ with this partition as diagonal, by choosing the off diagonal entries so that the 
vector states in $F$ are close to the trace (so close to being $0$ on $e_{ij}$ whenever $i\neq j$). All this is done recursively, in the style of [P18b]. The resulting 
inductive limit of matrix units gives rise to a UHF algebra $R_0\subset M$ with the property that all unit vectors 
in $\Cal H$ asymptotically  implement a trace on it. This readily  implies that its weak closure $R=\overline{R_0}^w$ gives a normal representation of the hyperfinite II$_1$ factor,  
while its diagonal $D$  is a MASA in $\Cal M$, showing that $R'\cap \Cal M=\Bbb C$. 

When applied to the particular case  $M=\Cal M$, the above theorem shows that, strikingly enough, the hyperfinite II$_1$ factor embeds ergodically 
into ANY continuous factor:    

 \proclaim{1.3. Corollary} Any separable continuous factor $\Cal M$ contains an ergodic copy of the hyperfinite $\text{\rm II}_1$ factor, 
 $R\hookrightarrow \Cal M$, and can be embedded ergodically into the unique AFD $\text{\rm II}_\infty$ factor, 
 $\Cal M \hookrightarrow R^\infty=R \overline{\otimes} \Cal B(\ell^2\Bbb N)$. 
 \endproclaim

The above corollary complements results in [P81a], which covered the case $\Cal M$ is II$_1$ or III$_\lambda$, $0<\lambda <1$, as well as results in 
[P83], [P84], [L84], which showed that if $\Cal M$ is III$_0$ or III$_1$, then it contains an irreducible AFD type III factor.  

The rest of the paper is organized as follows. In Section 2 we formally define the various versions of ergodicity 
and prove some basic properties. In Section 3 we prove Theorem 1.1. 
Then in Section 4 we prove some  criteria for an increasing sequence of finite partitions $A_n$ in a von Neumann 
algebra $\Cal M$ to generate a MASA $A=\overline{\cup_n A_n}^w$. In Section 5 we discuss the flatness property for MASAs 
while in Section 6 we prove Theorem 1.2.  In Section 7 we make some comments and state several open problems. 

While the results in this paper concern general von Neumann algebras, we will only need basic facts in this area, 
which can be found in  the Murray-von Neumann initial papers ([MvN36], [MvN43]) and the early monographs ([D57], [S71]). 
In particular, we will not use the Tomita-Takesaki theory in our arguments, and  
results on the structure and classification of type III factors in ([C73], [CT77], [C85], [H87])  
will only be used as ``black-boxes''. All of this background material  can be found in the   comprehensive  monograph ([T03]).  For basics in II$_1$ factors, 
notably on the hyperfinite II$_1$ factor, we refer the reader to ([AP17]).

\heading 2. Ergodicity properties for inclusions of factors
\endheading

\noindent
{\bf 2.1. Definition.} An inclusion of von Neumann algebras $M \subset \Cal M$ has the {\it weak relative Dixmier property}, or the {\it MV-property},   
if $\overline{\text{\rm co}}^w\{uxu^* \mid u\in \Cal U(M)\}\cap M'\cap \Cal M \neq \emptyset$, $\forall x\in \Cal M$. 
(N.B.: the abbreviation MV comes from ``mean value''.) 

\vskip.05in
\noindent
{\bf 2.2. Definitions}. Let $M \subset \Cal M$ be an inclusion  factors. 

$1^\circ$ $\Cal N\subset \Cal M$    is {\it irreducible}, or is {\it ergodic},  
if it satisfies the trivial relative commutant condition $\Cal N'\cap \Cal M=\Bbb C$.  

$2^\circ$ $M\subset \Cal M$ is {\it MV-ergodic} if $\overline{\text{\rm co}}^w\{uxu^* \mid u\in \Cal U(M)\}\cap \Bbb C1\neq \emptyset$, $\forall x\in \Cal M$. 
Also, If $\Cal U_0\subset \Cal U(M)$ is a subgroup and $\overline{\text{\rm co}}^w\{uxu^* \mid u\in \Cal U_0\}\cap \Bbb C1\neq \emptyset$, $\forall x\in \Cal M$, 
then we'll say that  $M\subset \Cal M$ is {\it MV-ergodic w.r.t. $\Cal U_0.$} If $M$ is II$_1$, then $M\subset \Cal M$ is {\it stably MV-ergodic} if $pMp \subset p\Cal Mp$ is MV-ergodic for any 
projection $p\in M$.

$3^\circ$ $M\subset \Cal M$ is {\it MASA-ergodic} if $M$ contains an abelian $^*$-subalgebra $A \subset M$ 
that's maximal abelian in $\Cal M$, i.e., $A'\cap \Cal M = A$. (N.B.: note that the condition $A'\cap \Cal M=A$ trivially implies $\Cal N \subset \Cal M$ is ergodic,  
because $\Cal N'\cap \Cal M \subset \Cal N' \cap (A'\cap \Cal M) =  \Cal N'\cap A \subset \Cal N'\cap \Cal N=\Bbb C$.) 

$4^\circ$ $M\subset \Cal M$ is $R$-{\it ergodic} if there exists a copy of the hyperfinite II$_1$ subfactor $R\subset M$ 
that's ergodic in $\Cal M$.  

Note that in case $M$ is of type II$_1$, then both MASA-ergodicity and $R$-ergodicity for embeddings $M \subset \Cal M$ are ``stable properties'', in that if $M\subset \Cal M$ has 
one of this property, then so does $pMp\subset p\Cal M p$, $\forall p\in \Cal P(M)$. 

\proclaim{2.3. Proposition} If an inclusion  of factors $M \subset \Cal M$ is either $R$-ergodic or MASA-ergodic, then 
it is MV-ergodic $($even stably MV-ergodic in case $M$ is $\text{\rm II}_1)$. If $M\subset \Cal M$ is MV-ergodic, then it is ergodic. 
\endproclaim
\noindent
{\it Proof}. The last implication is trivial, while the first two are observations made in [H87]: 
If for an element $x\in \Cal M$ and 
a von Neumann subalgebra $B\subset M$ we denote 
$K^w_B(x):=\overline{\text{\rm co}}^w \{uxu^* \mid u\in \Cal U(B)\} \cap M$, then Dixmier's averaging theorem applied to the factor $M$ shows 
that any $y\in K^w_B(x)$ satisfies $K_M^n(y):=\overline{\text{\rm co}}^n \{vyv^* \mid v \in \Cal U(M)\}\cap \Bbb C\neq \emptyset$. But $K^n_M(y) 
\subset K^w_M(x)$, so if $B\subset M$ is such that $K^w_B(x)\neq \emptyset$ (which is the case for $B=R$ ergodic in $\Cal M$, 
or for $B$ abelian and a MASA in $\Cal M$), then we have $K_M^w(x)\neq \emptyset$. 

If in addition $M$ is II$_1$ and one considers the inclusion $pMp\subset p\Cal Mp$ for some 
projection $p\in M$, then one can take $p\in B$ in both the case $B\simeq R$ ergodic in $\Cal M$, or $B$ abelian and a MASA in $\Cal M$, 
and the first part applies.  
\hfill $\square$

\vskip.05in

We recall below J. Schwartz's observation in [Sc63] that if $K_M^w(x):=\overline{\text{\rm co}}^w \{uxu^* \mid u\in \Cal U(M)\}\cap M'\cap \Cal M\neq \emptyset$, 
$\forall x\in \Cal M$, then one can choose the elements in $K^w_M(x)$ to depend linearly on $x$ (see [Sc63] for the proof): 

\proclaim{2.4. Lemma} Let $M \subset \Cal M$ be an inclusion of von Neumann algebras and $\Cal U_0\subset \Cal U(M)$ a 
subgroup  satisfying $\Cal U_0''=M$. The following conditions are equivalent: 

\vskip.05in 
$1^\circ$ Given any $x\in \Cal M$, we have $\overline{\text{\rm co}}^w \{uxu^* \mid u\in \Cal U_0\}\cap (M'\cap \Cal M)\neq \emptyset$.  

\vskip.05in 
$2^\circ$ For any finite $n$-tuple $(x_j)_j \subset (\Cal M)_1$, any finite set $F\subset \Cal M_*$ and any $\varepsilon >0$, 
there exist finitely many unitary elements $u_1, ..., u_m \in \Cal U_0$ and an $n$-tuple 
$(y_j)_j \in (M'\cap \Cal M)_1$ such that $|\frac{1}{m}\sum_{i=1}^m \varphi(u_i (x_j - y_j)u_i^*)| \leq \varepsilon$, $1\leq j \leq m$, $\varphi \in F$. 

\vskip.05in 
$3^\circ$ There exists a norm one projection 
of $\Cal M$ onto $M'\cap \Cal M$ that's in the weak closure of the convex set of c.p. maps $T_U:\Cal M \rightarrow \Cal M$, where 
$T_U(y)=\sum_i \alpha_i u_iyu_i^*$, $y\in \Cal M$, with $U=(u_i)_i\subset \Cal U(M)$ and $\alpha_i$ positive scalars summing up to $1$.

\endproclaim

\proclaim{2.5. Proposition} Let  $M \subset \Cal M$ be an inclusion of von Neumann algebras and 
$\Cal U_0\subset \Cal U(M)$ a 
subgroup  with the property that $\Cal U_0''=M$. The following conditions are equivalent: 

\vskip.05in 

$1^\circ$ Given any $x\in \Cal M$, we have $\overline{\text{\rm co}}^w \{uxu^* \mid u\in \Cal U_0\}\cap (M'\cap \Cal M)\neq \emptyset$. 

\vskip.05in

$2^\circ$ Given any finite set $F\subset \Cal M_*$, with $\varphi_{|M'\cap \Cal M}=0$, $\forall \varphi \in F$, 
the norm closure in the Banach space $(\Cal M_*)^{\oplus F}$ 
of the convex hull of the set $\{(u \cdot \varphi \cdot u^*)_{\varphi\in F} \mid u \in \Cal U_0 \}$, contains $0=(0)_{\varphi \in F}$. 

\vskip.05in 

$3^\circ$ Given any finite set $F\subset \Cal M_*$, with $\varphi_{|M'\cap \Cal M}=0$, $\forall \varphi\in F$, and any $\varepsilon >0$, 
there exist unitary elements $u_1, ..., u_n \in \Cal U(M)$ such that $\|\frac{1}{n} \sum_i u_i\cdot \varphi \cdot u_i^*\|\leq \varepsilon$, 
$\forall \varphi \in F$.  

\vskip.05in
Moreover, if these conditions are satisfied, then they are also satisfied by any other group $\Cal U_1\subset \Cal U(M)$ 
that has the same closure in the $s^*$-topology as $\Cal U_0$. 

\endproclaim
\noindent
{\it Proof}. One clearly has $2^\circ \Leftrightarrow 3^\circ$. 

To prove $1^\circ \Rightarrow 3^\circ$, assume 
$M\subset \Cal M$ satisfies MV property but  that there exists $F\subset \Cal M_*$ finite, with $\varphi(M'\cap \Cal M)=0$, $\forall \varphi \in F$, such that $\Cal K_F
:=\overline{\text{\rm co}}^n \{ (u\cdot \varphi \cdot u^*)_{\varphi \in F} \mid u\in \Cal U_0 \} \subset (\Cal M_*)^{\oplus F}$ does not contain 
the $F$-tuple $0$. 

Since the dual of the Banach space  $(\Cal M_*)^{n}$ is $\Cal M^{n}$, with the duality 
given by $\langle (\psi_j), (x_j)_j \rangle =\sum_j \psi_j(x_j)$, $\forall (\psi_j)_j \in (\Cal M_*)^{\oplus F}$, 
$(x_j)_j \in \Cal M^{\oplus F}$, by the Hahn-Banach separation theorem there exists $(x_\varphi)_{\varphi \in F}\in (\Cal M^{\oplus F})_1$ 
that separates the closed convex set $\Cal K_F \subset (\Cal M_*)^{\oplus F}$ from the compact set $\{0\}$, i.e., $\exists c>0$ 
such that 

$$
\text{\rm Re} \sum_{\varphi \in F}\varphi (ux_\varphi u^*) \geq c, \forall u\in \Cal U_0.
$$ 
But by Lemma 2.4,  
there exist $u_1, ..., u_m \in \Cal U_0$ and $(y_\varphi)_{\varphi \in F}  \subset (M'\cap \Cal M)_1$, such that 
$|\varphi(\frac{1}{m}\sum_{j=1}^m (u_j (x_\varphi - y_\varphi) u_j^*)| \leq c/2|F|$, $\forall \varphi \in F$. Since $\varphi(y_\varphi)=0$, 
when combined with the above inequality applied to $u=u_j$, $1\leq j \leq m$, this gives: 
$$
c/2 \geq \sum_{\varphi \in F} |\varphi(\frac{1}{m}\sum_{j=1}^m u_j x_\varphi  u_j^*)|
\geq \text{\rm Re} \sum_{\varphi \in F} (\frac{1}{m}\sum_{j=1}^m \varphi (u_jx_\varphi u_j^*) \geq c,
$$
a contradiction.

$3^\circ \Rightarrow 1^\circ$. If $M \subset \Cal M$ does not satisfy $1^\circ$  then there exists $x\in (\Cal M)_1$ 
such that $\Cal K_x:=\overline{\text{\rm co}}^{\sigma(\Cal M, \Cal M_*)} \{uxu^* \mid u \in \Cal U_0\}$ 
satisfies $\Cal K_x \cap (M'\cap \Cal M) = \emptyset$. 
By the Hahn-Banach theorem, since $\Cal K_x$ is $\sigma(\Cal M, \Cal M_*)$-compact and $M'\cap \Cal M$ is weakly closed, 
there exists a functional $\varphi\in \Cal M_*$ such that $\varphi$ is equal to zero 
on the weakly closed space $M'\cap \Cal M$ and $\text{\rm Re}\varphi(y) \geq c $, for some $c>0$ and all 
$y\in \Cal K_x$. In particular, given any $u_1, ..., u_m \in \Cal U_0$, we have 
$$
\|\frac{1}{m}\sum_{j=1}^m u_j\cdot \varphi \cdot u_j^*\| \geq |\frac{1}{m}\sum_{j=1}^m \varphi(u_jxu_j^*)| \geq 
\text{\rm Re}(\frac{1}{m}\sum_{j=1}^m \varphi(u_jxu_j^*)) \geq c. 
$$
But by $3^\circ$, one can take the unitaries $u_j\in\Cal U_0$ so  that the left hand term of these inequalities is arbitrarily small, 
say less than $c/2$, giving us a contradiction. 

\hfill $\square$

\proclaim{2.6. Corollary} Let  $M \subset \Cal M$ be an inclusion of factors acting on a Hilbert space $\Cal H$. The following conditions are equivalent: 

\vskip.05in 

$1^\circ$ $M \subset \Cal M$ is MV-ergodic. 

\vskip.05in

$2^\circ$ Given any finite set $F\subset \Cal M_*$, with $\varphi(1)=0$, $\forall \varphi\in F$, and any $\varepsilon >0$, 
there exist unitary elements $u_1, ..., u_m \in \Cal U(M)$ such that $\|\frac{1}{m} \sum_{i=1}^m u_i\cdot \varphi \cdot u_i^*\|\leq \varepsilon$, 
$\forall \varphi \in F$.  

\vskip.05in

$3^\circ$ Given any finite set of unit vectors $X\subset \Cal H$ and any $\varepsilon >0$, 
there exist unitary elements $u_1, ..., u_n \in \Cal U(M)$ such that $\|\frac{1}{n} \sum_i u_i\cdot (\omega_\xi - \omega_\zeta)_{|\Cal M} \cdot u_i^*\|\leq \varepsilon$, 
$\forall \xi, \zeta \in X$, where for a unit vector $\xi\in \Cal H$, $\omega_\xi$ denotes the state $\langle \ \cdot \ \xi, \xi \rangle$ implemented by $\xi$. 
\endproclaim
\noindent
{\it Proof}. Conditions $1^\circ$, $2^\circ$ are just $1^\circ$ and respectively $3^\circ$ in Proposition 2.5, in the case $M'\cap \Cal M=\Bbb C$, showing that 
$1^\circ \Leftrightarrow 2^\circ$. 

One clearly has $2^\circ \Rightarrow 3^\circ$. Then $3^\circ \Rightarrow 2^\circ$ follows trivially 
 from the fact that any functional $\varphi\in \Cal M_*$ with $\varphi(1)=0$ has 
real and imaginary parts also vanishing at $1$, and each selfadjoint functional that vanishes at $1$ is a scalar multiple of the difference between normal states on $M$, which  
are restrictions of normal states on $\Cal B(\Cal H)$, which in turn are norm limits of convex combination of vector states on $\Cal B(\Cal H)$, 
showing that in order to verify $2^\circ$, it is sufficient to check it for functionals that are restrictions to $\Cal M$ of differences of 
vector states in $\Cal B(\Cal H)$ 
\hfill $\square$

\vskip.05in
\noindent
{\bf 2.7. Remark.} Note that, due to Proposition 2.5, if $M\subset \Cal M$ is MV-ergodic, then it is MV-ergodic 
with respect to any subgroup $\Cal U_0\subset \Cal U(M)$ that's $s^*$-dense in $\Cal U(M)$. However, MV-ergodicity may fail to hold true 
with respect to subgroups $\Cal U_0$ that only satisfy  
$\Cal U_0''=M$, as the following  example shows: let $\Gamma$ be a non-amenable group and 
$\Gamma \curvearrowright R=(\Bbb M_{2}(\Bbb C), tr)^{\overline{\otimes} \Gamma}$  the Bernoulli $\Gamma$-action with base 
$(\Bbb M_{2}(\Bbb C), tr)$; let $P=R \rtimes \Gamma$ with $\{u_g\}_g \subset P$ the canonical unitaries implementing the 
$\Gamma$-action on $R$; let $M=\Cal M=\Cal B(L^2P)$ and define  
$\Cal U_0$ to be the subgroup of unitary elements on $L^2M$ 
generated by $\{u_g\}_g$ and the left and right multiplication by unitaries of $R$; then $M\subset \Cal M$ is MV-ergodic 
because of Dixmier's averaging theorem applied to the factor $M=\Cal M=\Cal B(L^2P)$ and we have $\Cal U_0''=M$, but the averaging by $\Cal U_0$ on elements in 
the abelian von Neumann algebra $\{u_g e_R u_g^*\}_g'' \simeq \ell^\infty(\Gamma)\subset \Cal B(L^2P)=\Cal M$ amounts 
to averaging by $\{u_g\}_g$; if there would be a state on $\ell^\infty\Gamma$ that's a weak limit of such averaging, then this would imply $\Gamma$ is amenable, a contradiction 
(see also the proof of Theorem 3.2 below).

\heading 3.  Two classes of counterexamples 
\endheading 

While $R$-ergodicity trivially implies MV-ergodicity (cf. 2.3), which in turn implies plain ergodicity, we'll show in this section 
that the opposite implications fail in general. 

Let us first notice that a result in [P18a] implies that whenever an ergodic II$_1 \subset \text{\rm II}_\infty$ inclusion of factors $M\subset \Cal M$ arises 
as the basic construction of an ergodic  quasi-regular inclusion, $N\subset M \subset \langle M, e_N \rangle = \Cal M$, with the $C^*$-tensor category generated by 
$_N L^2M_N$ amenable (see 2.12 in [P18a] for detailed definitions), then $M\subset \Cal M$ does follow $R$-ergodic. 

\proclaim{3.1. Proposition} Let $N \subset M$ be an ergodic quasi-regular inclusion of $\text{\rm II}_1$ factors. Assume  the 
irreducible Hilbert-bimodules contained in $_NL^2(M)_N$ generate an amenable concrete $C^*$-tensor category $\Cal G$, in the sense of $2.12$ in $\text{\rm  [P18a]}$. 
If $N\subset M \subset \Cal M=\langle M, e_N \rangle$ is the associated basic construction, then $M \subset \Cal M$ is $R$-ergodic $($and thus 
stably MV-ergodic as well$)$. 
\endproclaim 
\noindent
{\it Proof}. By (Theorem 2.12 in [P18a]), there exists a hyperfinite subfactor $Q\subset N$ that's ``normalized'' by $\Cal G$ and satisfies $Q'\cap \Cal M=N'\cap \Cal M$. 
This amounts to the fact that there exists a hyperfinite subfactor $R\subset M$ that contains $Q$ such that the inclusion  of $Q\subset R$ into $N\subset M$ is a non-degenerate 
commuting square, with $e_N$ implementing the basic construction $Q \subset R \subset \langle R, Q \rangle=\langle R, e_N\rangle$, and 
with $Q'\cap \Cal M=N'\cap \Cal M$, $R'\cap \Cal M= M'\cap \Cal M=\Bbb C$. Thus, $M\subset \Cal M$ is $R$-ergodic. 
\hfill $\square$

\proclaim{3.2. Theorem} Let $N$ be a $\text{\rm II}_1$ factor, $\Gamma \curvearrowright N$ a free action of a discrete group, $M=N \rtimes \Gamma$ 
its crossed product factor and $N \subset M \subset \Cal M= \langle M, e_N \rangle$ the associated basic construction. The following conditions are equivalent: 

\vskip.05in

$1^\circ$ $M \subset \Cal M$ is MV-ergodic

$2^\circ$ $M\subset \Cal M$ is $R$-ergodic. 

$3^\circ$ The group $\Gamma$ is amenable. 

$4^\circ$ The regular inclusion $N \subset M$ is co-amenable. 
\endproclaim 
\noindent
{\it Proof}. $3^\circ \Rightarrow 1^\circ$ The inclusion $N \subset \Cal M$  always has the MV property (see e.g. 1.4 in [P98]). Let $x\in \Cal M$ and 
$x_0 \in \overline{\text{\rm co}}^w \{uxu^* \mid u\in \Cal U(N)\} \cap (N'\cap \Cal M)$. 
Note that $N'\cap \Cal M= \ell^\infty(\Gamma)$ with the canonical unitaries $\{u_g\}_g\subset M$ implementing 
the left translation on $\ell^\infty \Gamma$ via the Ad-action.  Thus,  if $\Gamma$ is amenable, then by averaging 
$x_0$ over the canonical unitaries with respect to the invariant mean on $\Gamma$ gives a scalar as weak limit. 

$1^\circ \Rightarrow 3^\circ$ If $M \subset \Cal M$ is MV-ergodic, then by using Lemma 2.4 it follows that there exists a state $\varphi$ on $\ell^\infty\Gamma$ that's 
obtained as a pointwise weak limit of averaging by unitaries in $ M$. Since all unitaries in $M$ commute with all $u_g^{op}$ and $\text{\rm Ad}(u_g^{op})$ normalize 
$\Cal M$, one has $\varphi = \varphi \circ \text{\rm Ad}(u_g^{op})$, $\forall g$. In particular, this implies 
$\varphi(f)=\varphi(u_g^{op}f{u_g^{op}}^*)$ for any $f \in \ell^\infty \Gamma=N'\cap \Cal M$ 
and $g\in \Gamma$. But $\text{\rm Ad}(u_g^{op})$ implements the right translation of $f$ by $g$, $f_g$. Thus, $\varphi(f)=\varphi(f_g)$, 
showing that $\varphi$ is a (right) invariant mean on $\Gamma$, i.e., $\Gamma$ is amenable.  

$3^\circ \Rightarrow 2^\circ$ is a a particular case of Proposition 3.1 while $2^\circ \Rightarrow 1^\circ$ is trivial (see Proposition 2.3). 

The equivalence $3^\circ \Leftrightarrow 4^\circ$ is well known (see e.g. [P01]).  

\hfill $\square$ 

\proclaim{3.3. Corollary} If $\Gamma$ is a nonamenable group, $\Gamma \curvearrowright R=\overline{\otimes}_{g\in \Gamma} (\Bbb M_2(\Bbb C), tr)_g$ 
is the Bernoulli $\Gamma$-action with base $(\Bbb M_2(\Bbb C), tr)$, and we let $R \subset M=R \rtimes \Gamma$, then 
the basic construction inclusion $M \subset \langle M, e_R \rangle \simeq R^\infty:= R\overline{\otimes} \Cal B(\ell^2\Bbb N)$ is ergodic but not MV-ergodic. 
\endproclaim 

We next combine results in [DPe19] and [GP96] to prove the existence of MV-ergodic inclusions of factors that 
are not $R$-ergodic. The examples we give are the consequence of an ``alternative'' between two possible inclusions, obtained by using the Das-Peterson 
notion of Poisson boundary of a II$_1$ factor [DPe19], which we recall here for the reader's convenience. 

Thus, let $M$ be a separable II$_1$ factor,  
$1\in \Cal U=\{u_n\}_n \subset M$ an at most countable self-adjoint set of unitary elements that generate $M$ as a von Neumann algebra, $\alpha=\{\alpha_n\}_n$ 
some positive ``weights'', satisfying $\sum_n \alpha_n=1$ and $\alpha_n=\alpha_m$ if $u_m^*=u_n$. Let 
$\varphi=\varphi_{\Cal U, \alpha}$ be the unital completely positive (ucp) map on $\Cal B(L^2M)$ given by $\varphi(T)=\sum_n \alpha_n u_n^{op} T {u_n^{op}}^*$, $T\in \Cal B(L^2M)$, 
and denote $\Cal B_\varphi$ its fixed points, $\{T\in \Cal B(L^2M)\mid \varphi(T)=T\}$. Note that $M\subset \Cal B_\varphi$. 
Endow this weakly closed self-adjoint subspace of $\Cal B(L^2M)$  
with a Choi-Effros product, given by $X \cdot Y:=P_\varphi(XY)$, $\forall T, S\in \Cal B_\varphi$, 
where $P_\varphi: \Cal B(L^2M)\rightarrow \Cal B(L^2M)$ is the ucp map  obtained by averaging over $\varphi$, i.e., 
$P_\varphi(T)=\underset{n \rightarrow \omega}\to{\lim} \frac{1}{n}\sum_{k=1}^n \varphi^k(T)$ ($\omega$ denotes here a free ultrafilter on $\Bbb N$). 

It is shown in [DPe19] that $\Cal B_\varphi$ with this product is a von Neumann factor, which is injective (or amenable), and thus 
AFD by Connes Theorem [C76],  that it is properly infinite whenever $M\neq R$, and that the inclusion of factors 
$M \subset B_\varphi$ has trivial relative commutant, i.e., it is ergodic.  We will call the inclusion of factors $M\subset \Cal B_\varphi$ 
the (Das-Peterson) {\it Poisson boundary of} $M$ corresponding to $\varphi$. 

An important {\it double ergodicity} result in [DPe19] shows that in fact $M \subset \Cal B_\varphi$ is MV-ergodic. More precisely, one has 
$\Cal B_\varphi \cap J\Cal B_\varphi J=\Bbb C$, where $J=J_M: L^2M \rightarrow L^2M$ is the canonical conjugation, $J(\xi)=\xi^*, \xi \in L^2M$. 
This amounts to the following: given any $T\in \Cal B_\varphi$, there exists an averaging of $T$ by $J \varphi J = \sum_n \alpha_n u_n \cdot u_n^*$ and its powers that converges weakly to a scalar.    
Or, with the above notations, $JP_\varphi J(T)\in \Bbb C1$, $\forall T\in \Cal B_\varphi$.

\proclaim{3.4. Theorem}  Let $M$ be an interpolated free group factor, i.e., $M=L(\Bbb F_t)$ for some $1< t \leq \infty$ $(\text{\rm [Dy93], [R92]})$. Let  
$M \subset \Cal B_\varphi$ be a Poisson boundary of $M$, corresponding to some $\varphi$ as above. 
Then either the Das-Peterson MV-ergodic inclusion $M \subset \Cal B_\varphi$ is not $R$-ergodic, or if it is $R$-ergodic and $L\subset M$ is a hyperfinite 
$\text{\rm II}_1$ subfactor such that $L'\cap \Cal B_\varphi=\Bbb C$, then $M\subset \langle M, e_L \rangle$ is MV-ergodic 
but not $R$-ergodic. 
\endproclaim
\noindent
{\it Proof}. Assume $M \subset \Cal B_\varphi$ is $R$-ergodic and that $L\subset M$ is a hyperfinite subfactor with $L'\cap \Cal B_\varphi=\Bbb C$. 
By conjugating with $J \cdot J$ and taking into account that $J L' J=\langle M, e_L \rangle$ and the definition of $J\Cal B_\varphi J$, this 
is equivalent to the fact that for any $T\in \langle M, e_L \rangle$, we have $JP_\varphi J(T)\in \Bbb C$, so an averaging by $J\varphi J=\sum_n \alpha_n u_n \cdot u_n^*$ 
and its powers can be taken to converge weakly to scalars. Thus, $M\subset \langle M, e_L \rangle$ is MV-ergodic. But if there exists some hyperfinite subfactor 
$R \subset M$ such that $R'\cap \langle M, e_L \rangle = \Bbb C$, equivalently $R\vee L^{op}=\Cal B(L^2M)$, then $[R \xi L]=L^2M$ for any $\xi \neq 0$ in $L^2M$, 
contradicting (Theorem 4.2 in [GP96]). 
\hfill $\square$

\heading 4. MASA criteria
\endheading 

In this section we prove several criteria for an abelian $^*$-subalgebra $A$ 
of a von Neumann algebra $\Cal N$ to be maximal abelian (i.e., a MASA) in $\Cal N$. They all 
start from a representation of $A$ as an increasing limit of finite dimensional subalgebras $A_n \subset A$ (which we will 
also call {\it finite partitions} of $1$ with projections in $A$) and will provide MASA conditions for $A\subset M$ in terms 
of ``local conditions'' for $A_n$. This will allow constructing MASAs $A\subset M$ recursively, from the finite data $A_n$, which are chosen so that to become  
``more and more maximal abelian'' in $\Cal N$, while at the same time to satisfy other conditions.

\proclaim{4.1. Lemma} Let $\Cal N$ be a separable von Neumann algebra and $A\subset \Cal  N$ an abelian 
von Neumann subalgebra. Let $A_n \subset A$ be an increasing sequence of finite dimensional von Neumann 
subalgebras that generate $A$, i.e., $\overline{\cup_n A_n}^w=A$. Then we have:

\vskip.05in

$1^\circ$ $\lim_n \|\varphi \circ E_{A_n'\cap \Cal N}\|=0$, for all $\varphi \in \Cal N_*$ with $\varphi(A'\cap \Cal N)=0$. 

\vskip.05in

$2^\circ$ $A$ is a MASA in $\Cal N$ if and only if 
$\lim_n \|\varphi \circ E_{A_n'\cap \Cal N}\|=0$, for all $\varphi \in \Cal N_*$ with $\varphi(A)=0$. 

\endproclaim
\noindent
{\it Proof}. $1^\circ$ Let $\varphi\in \Cal N_*$ with $\varphi(A'\cap \Cal N)=0$. Since $\Cal U_0=\cup_n \Cal U(A_n)$ is $s^*$-dense in $\Cal U(A)$, by condition $3^\circ$ 
in Proposition 2.5 
it follows that for any $\varepsilon >0$ there exist an $n_0$ and $u_1, ..., u_m \in \Cal U(A_{n_0}) \subset \Cal U_0$ 
such that $\| \frac{1}{m}\sum_{i=1}^m u_i \cdot \varphi \cdot u_i^*\| < \varepsilon$. But if $n \geq n_0$ then for any unitary element $u\in A_n$ 
(so in particular for all $u_i$) we have $(u \cdot \varphi \cdot u^*) \circ E_{A_n'\cap \Cal N} = \varphi \circ E_{A_n'\cap \Cal N}$. 
Thus, for all $n\geq n_0$ we have 
$$   \|\varphi \circ E_{A_n'\cap \Cal N}\| = \|(\frac{1}{m}\sum_{i=1}^m u_i \cdot \varphi \cdot u_i^*) \circ E_{A_n'\cap \Cal N}\| 
\leq \|\frac{1}{m}\sum_{i=1}^m u_i \cdot \varphi \cdot u_i^*\|\leq \varepsilon,  
$$
showing that $\lim_n \|\varphi \circ E_{A_n'\cap \Cal N}\|=0$.

\vskip.05in

$2^\circ$ Since $A$ is a MASA iff $A'\cap \Cal N = A$, this part follows immediately from $1^\circ$. 
\hfill $\square$

\proclaim{4.2. Lemma} Let $\Cal N\subset \Cal B(\Cal H)$ be a von Neumann algebra acting on a separable Hilbert space, $A_n\subset \Cal N$ an increasing 
sequence of finite dimensional abelian von Neumann subalgebras and $A=\overline{\cup_n A_n}^w$. 
\vskip.05in
$1^\circ$ Given any $\xi\in \Cal H$, we have $[A_n\xi]\nearrow [A\xi] \in A'$. 

$2^\circ$ Let $\{e^n_j\}_j$ denote the minimal projections in $A_n$. Given any $x\in (A \vee \Cal N')_1$, there exist $x_j^n\in (\Cal N')_1$ 
such that $x$ is the limit in the strong operator topology of the sequence $\{\sum_j e^n_j x^n_j\}_n$. Moreover, if $x\geq 0$ then one can take $x_j^n$ positive. 
\endproclaim
\noindent
{\it Proof}. Part $1^\circ$ is trivial and part $2^\circ$ is an immediate consequence of the fact that the $^*$-algebra $\cup_n A_n \vee \Cal N'$ is weakly dense 
in $A \vee \Cal N'$ and Kaplanski's density theorem.  

\hfill $\square$

\proclaim{4.3. Lemma} Let $\Cal N\subset \Cal B(\Cal H)$ be a von Neumann algebra  represented on a separable Hilbert space $\Cal H$, 
with $\xi\in \Cal H$ a  separating unit vector. Let $\{\eta_m\}_m$ be a total sequence in $\Cal H$. 
Let $A\subset \Cal N$ be an abelian von Neumann subalgebra. The following conditions are equivalent: 

\vskip.05in

$1^\circ$ $A$ is maximal abelian in $\Cal N$, i.e. $A'\cap \Cal N=A$. 

\vskip.05in 

$2^\circ$ There exist an increasing sequence of 
finite partitions $\{e^n_j\}_{1\leq j \leq k_n} \subset \Cal P(A)$, $n\geq 1$, and positive elements $\{x^n_j\}_{1\leq j \leq k_n}\subset (\Cal N')_1$, $n\geq 1$, 
such that if we denote $A_n=\sum_j \Bbb C e_j^n$ then $\lim_n \|[A_n \xi](\eta_m)-\sum_{j=1}^{k_n} e^n_j x^n_j ( \eta_m)\|=0$, $\forall m$. 

\vskip.05in 

$3^\circ$ Given any increasing sequence of 
finite partitions $\{e^n_j\}_{1\leq j \leq k_n}$ $\subset \Cal P(A)$, $n\geq 1$, with $A_n=\sum_j \Bbb Ce^n_j$ satisfying $\overline{\cup_n A_n}^w=A$, 
there exist positive elements $\{x^n_j\}_{ 1\leq j \leq k_n}$ $ \subset (\Cal N')_1$, $n\geq 1$, 
such that  $\underset{n \rightarrow \infty}\to{\lim} \|[A \xi](\eta)-\sum_{j=1}^{k_n} e^n_j x^n_j (\eta)\|=0$, $\forall \eta \in \Cal H$.

\endproclaim
\noindent
{\it Proof}.  $1^\circ \Rightarrow 3^\circ$ If $A$ is a MASA in $\Cal N$, then $A'\cap \Cal N = A$. By taking commutants,  this implies 
$A \vee \Cal N' = A'$. Since $[A\xi] \in A'$, it follows that $[A\xi] \in A \vee \Cal N'$. By Lemma 4.2.2$^\circ$,  there exist positive elements  
$x^n_j \in (\Cal N')_1$ such that the sequence $\{\sum_j e^n_j x^n_j\}_n$ converges  strongly to $[A\xi]$. 

$3^\circ \Rightarrow 2^\circ$ is trivial, by taking into account 4.2.1$^\circ$. 

Finally, if we assume $2^\circ$, then by Lemma 4.2 it follows that the projection $[A\xi]$, 
which belongs to the commutant of $A$ in $\Cal B(\Cal H)$, lies in fact in $A \vee \Cal N'$.  Let $\tilde{A}\subset \Cal N$ be a MASA that contains 
$A$. Then we have $[\tilde{A}\xi]\in \tilde{A}'=\tilde{A}\vee \Cal N' \supset A \vee \Cal N'$. Thus, both $[A\xi], [\tilde{A}\xi]$ belong to $\tilde{A}\vee \Cal N'$. 
Since any element in the set $[\tilde{A}\xi](\tilde{A}\vee \Cal N')[\tilde{A}\xi]$ commutes with $\tilde{A}$, we also have   $[\tilde{A}\xi](\tilde{A}\vee \Cal N')[\tilde{A}\xi] 
=\tilde{A}[\tilde{A}\xi]$. It thus follows that $[A\xi]$ belongs to $\tilde{A}[\tilde{A}\xi]$, which is an abelian von Neumann algebra with cyclic and separating vector $\xi$ 
on the Hilbert space $\overline{\tilde{A}\xi}$. This Hilbert space coincides with the standard representation of $\tilde{A}$ corresponding to the faithful normal state $\tau$ implemented by the vector $\xi$. 
But then $A\subset \tilde{A}$ and $[A\xi]=[\tilde{A}\xi]$ forces $A=\tilde{A}$ (see e.g., [AP17]), showing that $2^\circ \Rightarrow 1^\circ$. 
\hfill $\square$

\proclaim{4.4. Lemma} Assume $A$ is a MASA in a von Neumann algebra $\Cal N \subset \Cal B(\Cal H)$ with a separating vector $\xi\in \Cal H$. 
Let $F\subset \Cal H$ be a finite set and $\delta >0$. 
Assume $A_0\subset A$ is a finite dimensional von Neumann subalgebra, with minimal projections $\{e_k\}_k$, and $\{x_k\}_k\subset (\Cal N')_1$  are such that 
$\|[A_0\xi](\eta) - [A\xi](\eta)\| + \|[A\xi](\eta) - \sum_j e_jx_j(\eta)\|< \delta$, $\forall \eta \in F$. If $B_0 \subset A$ is a finite dimensional von Neumann subalgebra that contains $A_0$ 
and $\{e_{k,j}\}_{k,j}$ are its minimal projections with $\sum_je_{k,j}=e_k$, $\forall k$, then 
$\|[B_0\xi](\eta)- [A\xi](\eta)\| + \|[A\xi](\eta) - \sum_{k,j} e_{k,j}x_k(\eta)\|< \delta$, $\forall \eta \in F$. 
\endproclaim
\noindent
{\it Proof.} Since $[A_0\xi]\leq [B_0\xi]\leq [A\xi]$, it follows that for all $\eta \in F$ we have 
that 
$$
\|[B_0\xi](\eta) - [A\xi](\eta)\| + \|[A\xi](\eta) - \sum_{k,j} e_{k,j}x_k(\eta)\| 
$$
$$
\leq \|[A_0\xi](\eta) - [A\xi](\eta)\| + \|[A\xi](\eta) - \sum_{k,j}e_{k,j}x_k(\eta)\| 
$$
$$
=  \|[A_0\xi](\eta) - [A\xi](\eta)\| + \|[A\xi](\eta) - \sum_{k}e_{k}x_k(\eta)\|< \delta. 
$$
\hfill $\square$

\proclaim{4.5. Lemma} Let $\Cal N\subset \Cal B(\Cal H)$ be a properly infinite factor represented on a separable Hilbert space $\Cal H$  
 and $A\subset \Cal N$ a MASA. 

\vskip.05in 

$1^\circ$ If the atomic part of $A$ is either $0$ or infinite dimensional, 
then there exists a diffuse von Neumann subalgebra $B\subset A$ with all its non-zero projections infinite $($and thus equivalent to $1)$  in $\Cal N$. 

\vskip.05in

$2^\circ$ Given any $n\geq 1$ there exists a partition of $1$ with projections $p_1, ..., p_n\in A$ such that all $p_i$ are infinite in $\Cal N$, 
and thus $p_i \sim p_j$, $\forall i,j$.

\endproclaim 
\noindent
{\it Proof}. $1^\circ$ Note first that if $\Cal N$ is of type III of II$_\infty$ then $A$ is diffuse, so its atomic part is $0$. 

If $\Cal N$ is of type III then we can just take $B=A$.  If $\Cal N$ is of type II$_\infty$, then let $f\in A$ be the supremum over all finite projection in $A$. If $f=0$ then we can 
just take $B=A$. If $f=1$, then we can take a partition of $1$ with projections $p_k\in A$ of trace $Tr_{\Cal N}(p_k)=1$, $\forall k\geq 1$, then choose isomorphisms $\theta_k: Ap_1\simeq Ap_k$ 
that preserve the trace,  and define $B=\{\sum_k \theta_k(a) \mid a\in Ap_1\}$. Finally, if both $f, 1-f\neq 0$, then take any identification of the (separable) diffuse abelian 
von Neumann algebras $\theta: Af\simeq A(1-f)$ and define $B=\{a + \theta(a) \mid a\in Af\}$. 

In case $\Cal N$ is of type I, then let  $f\in \Cal P(A)$ be the projection with $Af$ atomic and $A(1-f)$ diffuse. If $f=1$, it means $A$ is atomic, and we can take any diffuse quotient $B\simeq L^{\infty}([0,1])$ 
of $A\simeq \ell^\infty\Bbb N$. If $f=0$, it means $A$ itself is diffuse and all its projections are infinite so we can just take $B=A$. If $f\neq 0,1$, then its atomic part $Af$ is infinite 
dimensional so it contains a diffuse von Neumann subalgebra $B_0\subset Af$, which is isomorphic with the diffuse von Neumann algebra $A(1-f)$. Letting $\theta:B_0 \simeq A(1-f)$ 
be such an isomorphism, we define $B=\{a + \theta(a) \mid a\in B_0\}$, which clearly satisfies the condition.

$2^\circ$ If the atomic part $Af$ of $A$ is either $0$ or infinite dimensional, then $1^\circ$ shows that one can take $B\subset A$ diffuse, and any partition of $1$ with $n$ non-zero 
projections in $B$ will satisfy the conditions. If $Af$ is non-zero but finite dimensional, then one can 
take any partition of $1-f$ with non-zero projections $p'_1,p_2, p_3,  ..., p_n$ in the diffuse part $A(1-f)$  and let $p_1=p'_1+f$. The partition $\{p_j\}_{j=1}^n \subset A$ 
will then have all its projections infinite (and thus equivalent) in $\Cal N$. 
\hfill $\square$ 

\proclaim{4.6. Lemma} Let $\Cal N\subset \Cal B(\Cal H)$ be a properly infinite factor represented on a separable Hilbert space $\Cal H$, 
with a separating unit vector  $\xi\in \Cal H$, and $A\subset \Cal N$ a diffuse MASA. Denote by $\phi$ the vector state implemented by $\xi$.

\vskip.05in 

$1^\circ$ If $\{e_j\}_{j=1}^n\subset \Cal P(A)$ is a finite partition, 
$F\subset \Cal H$ is a finite set and $\alpha > 0$, then 
there exists a finite partition $\{f_j\}_{j=1}^n\subset A$ such that: $(a)$ $f_i$ is infinite in $\Cal N$, $1\leq i\leq n$; 
$(b)$ $\|(e_i-f_i)(\eta)\|< \alpha$, $\forall \eta\in F$, $1\leq i \leq n$; $(c)$ $\phi(f_i)=\frac{k_i}{2^t}$, for some integers $t\geq 1$, $k_i\geq 1$, $1\leq i\leq n$; 
$(d)$ for all but possibly one $j\in \{1, ..., n\}$ we have $f_j\geq e_j$. 

\vskip.05in 

$2^\circ$  If $\{f_i\}_i \subset \Cal P(A)$ is a partition of $1$ with infinite projections such that $\phi(f_i)=\frac{k_i}{2^t}$, for some integers $t\geq 1$, $k_i\geq 1$, $1\leq i\leq n$,  
then there exists a refinement of $\{f_i\}_i$ with $2^t$-many projections $\{f^0_j\}_j\subset \Cal P(A)$ that are infinite in $\Cal N$ and satisfy $\phi(f^0_j)=2^{-p}$, $\forall j$. 
\endproclaim 
\noindent 
{\it Proof}. $1^\circ$ Since $\Cal N$ is infinite, at least one of the projections $e_1, e_2, ..., e_n$ is infinite, and we can assume $e_1$ is infinite. 
By the part $1^\circ$, there exists a diffuse von Neumann subalgebra $B\subset Ae_1$ with all its projections infinite. Since $\xi$ is separating, there exists $c>0$ such that 
if $g\in \Cal P(B)$ satisfies $\phi(g)<c$ then $\|g(\eta)\|< \alpha$, $\forall \eta \in F$. On the other hand, $(B, \phi)$ can be viewed as $(L^\infty([0,s]), \mu)$ where $\mu$ is 
the Lebesgue measure and $s=\phi(e_1)$. But then there clearly exist $n-1$ disjoint intervals $I_2, ..., I_n$ in $[0,s]$ of length $\mu(I_k)=s_k$ satisfying $\sum_k s_k<t$, with $s_k<c$ and  
$s_k+\phi(e_k)$ dyadic for each $k=2, 3, ..., n$. Letting $p_k$ be the projection in $B$ corresponding to the characteristic function of $I_k$, $f_k=e_k+p_k$, 
for $2\leq k \leq n$, and $f_1=1-\sum_{k=2}^n f_k$, we get a partition$\{f_k\}_{k=1}^n$ that clearly satisfies all the desired conditions. 

$2^\circ$ By Lemma 4.5, for each $i$ there exists a diffuse von Neumann subalgebra $B_i\subset Af_i$ with all projections infinite in $\Cal N$. Since $(B_i, \varphi_{|B_i})$ is 
isomorphic to $(L^\infty([0,\frac{k_i}{2^t}]), \mu)$, we can split each $f_i=1_{B_i}$ into a partition with $k_i$ projections in $B_i$ of $\phi$-measure equal to $2^{-t}$. 
\hfill $\square$
 
\vskip.05in
\noindent
{\bf 4.7. Remark.} When constructing recursively a MASA $A$ through approximation by finer and finer partitions $A_n \nearrow A$, as in Lemma 4.3, one often needs to argue 
that once a certain ``level of approximation'' $A_n$ is reached, any finer partition $A_n\subset A_n^0\subset A$ still satisfies that same degree of approximation. 
However, unlike the II$_1$ case, where this is indeed the case due to the existence of trace preserving expectations on all von Neumann subalgebras, in general 
one does not necessarily  have such good behavior. The more complicated local approximation in Lemma 4.4 circumvents this shortcoming, allowing the passage to 
arbitrary finer partitions.

\heading 5. Flatness of MASAs in  factors
\endheading 

In this section we prove a key property of MASAs $A$ in continuous factors $M$, that we call flatness, 
which requires that any finite set of normal states in $M$ can be simultaneously rotated by a unitary in $M$ so that when restricted to $A$ 
they implement ``almost the same measure''. By contrast, we'll show that MASAs in atomic  factors do not have this property, 
due to  lack of ``enough room''. 

\vskip.05in

\noindent
{\bf 5.1. Definition}.  Let $M$ be a  von Neumann factor and $A\subset M$ a MASA. We say that $A$ is {\it flat} in $M$ 
if for any finite set $F$ in the space of normal states on $M$, S$_{\text{\rm n}}(M)$, and any $\varepsilon > 0$, there exists a unitary element $u\in M$ such 
that $\|(\varphi - \psi)_{|uAu^*}\|\leq \varepsilon$, for any $\varphi, \psi \in F$ (equivalently, $\|(\varphi \circ \text{\rm Ad}(u) - \psi \circ \text{\rm Ad}(u))_{|A}\|\leq \varepsilon$,  $\forall \varphi, \psi \in F$). 

Note right away that MASAs in a finite dimensional factor $M\neq \Bbb C$ cannot have this property, because one can take $F$ 
a large finite set that's ``almost-dense'' in S$_{\text{\rm n}}(M)$, for which the above condition will of course not hold. 

\proclaim{5.2. Lemma} Let $M$ be an infinite dimensional von Neumann  factor and  $A\subset M$ a MASA. The following conditions are equivalent: 

\vskip.05in 
$1^\circ$ $A$ is flat in $M$; 

\vskip.05in

$2^\circ$ For any finite set $F_0 \subset \{\varphi \in M_*\mid \varphi(1)=0 \}$ and any $\varepsilon >0$ there exists a unitary element $u\in M$ such that 
$\|\varphi_{|uAu^*}\|< \varepsilon$, $\forall \varphi\in F$. 

\vskip.05in

$3^\circ$ There exists a subset of normal functionals $\Cal L \subset M_*$ vanishing at $1$ and satisfying $\overline{\text{\rm sp}}^n \Cal L = \{\varphi \in M_*\mid \varphi(1)=0 \}$ 
with the property that for any $ E_0 \subset \Cal L$ finite and any $\delta >0$,  there exists $ u\in \Cal U(M)$ such that 
$\|\psi_{|uAu^*}\|< \delta$, $\forall \psi \in E_0$.

\endproclaim
\noindent
{\it Proof}. $2^\circ \Rightarrow 3^\circ$ is clear, by just taking $\Cal L=\{\varphi\in M_* \mid \varphi(1)=0\}$. 

$3^\circ \Rightarrow 2^\circ$ Assume  $\Cal L\subset\{\varphi\in M_* \mid \varphi(1)=0\}$  satisfies condition $3^\circ$. 
Let $F_0=\{\varphi_i\}_i \subset \{\varphi\in M_* \mid \varphi(1)=0\}$ be a finite set and $\varepsilon >0$. Then 
 there exist $E_0 = \{\psi_j\}_j \subset \Cal L$ finite and scalars $c_{ij}$ such that $\|\varphi_i - \sum_j c_{ij}\psi_j\|< \varepsilon/2$, $\forall i$. 
By property $3^\circ$, given any $\delta >0$ there exists $u\in \Cal U(M)$ such that $\|c_{ij} {\psi_j}_{|uAu^*}\|< \delta=\varepsilon/2|E_0|$, $\forall i,j$. 
Thus, we have 
$$
\|{\varphi_i}_{|uAu^*}\| \leq \|\varphi_i - \sum_j c_{ij}\psi_j \| + \sum_j \|c_{ij} {\psi_j}_{|uAu^*}\| \leq \varepsilon/2 + \varepsilon/2=\varepsilon
$$

$2^\circ \Rightarrow 1^\circ$ is trivial, then  $1^\circ \Rightarrow 3^\circ$ follows from the fact that 
the set $\Cal L_0\subset M_*$ of normal functionals of the form $\varphi - \psi$, with $\varphi, \psi$ normal states on $M$,  
has norm-dense linear span in  $\{\varphi \in M_*\mid \varphi(1)=0 \}$. 

\hfill $\square$ 

\proclaim{5.3. Lemma} Let $M$ be a  factor, $A\subset M$ a MASA and $N\subset M$ a subfactor that contains $A$. 
If $A$ is flat in $N$ then it is flat in $M$. 
\endproclaim
\noindent
{\it Proof}. This is trivial, by the definition. 

\hfill $\square$

\proclaim{5.4. Lemma} Let $M$ be an infinite dimensional von Neumann factor. The following conditions are equivalent: 

\vskip.05in

$1^\circ$  $M$ has a flat MASA

\vskip.05in 

$2^\circ$ Any MASA in $M$ is flat. 

\vskip.05in 

$3^\circ$ Given any finite set $F_0\subset \{\varphi \in M_*\mid \varphi(1)=0\}$ and any $\varepsilon >0$, there exists a finite partition of $1$ with projections 
$\{p_j\}_j\subset \Cal P(M)$ such that $p_i \sim p_j$, $\forall i,j$, and $\|\sum_i p_i \cdot \varphi \cdot p_i \|< \varepsilon$, $\forall \varphi \in F_0$. 
\endproclaim
\noindent
{\it Proof}. $3^\circ \Rightarrow 2^\circ$. Let $A\subset M$ be a MASA. Let $F_0\subset \{\varphi \in M_*\mid \varphi(1)=0\}$,  $\varepsilon >0$.   
By  condition $3^\circ$ we can find a finite partition of $1$ with mutually equivalent projections $\{p_j\}_j \subset M$ such that 
$\|\sum_j p_j \cdot \varphi \cdot p_j \|< \varepsilon$, $\forall \varphi \in F_0$. Since $A\subset M$ is a MASA in an infinite dimensional factor, by Lemma 4.5.2$^\circ$ 
it follows that $A$ contains a partition of $1$ with infinite, mutually equivalent 
projections $\{p'_j\}_j\subset A$ (same number as the partition $\{p_j\}_j$). But then there exists a unitary element $u\in M$ such 
that $up'_ju^*=p_j$, $\forall j$. On $uAu^*$ we then have for any $\varphi \in F_0$ 
$$
\|\varphi_{|uAu^*}\|\leq \|\varphi_{|\sum_i p_iMp_i}\| = \|\sum_i p_i \cdot \varphi \cdot p_i\| < \varepsilon.
$$
By Lemma 5.2, this shows that $A$ is flat in $M$. 

One clearly has $2^\circ \Rightarrow 1^\circ$. 

$1^\circ \Rightarrow 3^\circ$. Let $F_0, \varepsilon$ be given as in $3^\circ$. 
Let $B\subset M$ be a flat MASA and $u\in \Cal U(M)$ be so that $\|\varphi_{|uBu^*}\|< \varepsilon/2$, $\forall \varphi \in F_0$.  
Then $A=uBu^*$ is a MASA with $\|\varphi_{|A}\|< \varepsilon/2$, $\forall \varphi \in F_0$. By Hahn-Banach, it follows that 
for each $\varphi_k\in F_0 = \{\varphi_k\}_k$ there exists $\varphi'_k \in M_*$ such that ${\varphi'_k}_{|A}={\varphi_k}_{|A}$ and $\|\varphi'_k\|< \varepsilon/2$. 

Thus, if we denote $\psi_k=\varphi_k-\varphi'_k$ then $\psi_k$ vanishes on  $A$ and satisfies $\|\varphi_k-\psi_k\| = \|\varphi'_k\| < \varepsilon/2$.  

Let now $A_n \nearrow A$ be an increasing sequence of finite dimensional subalgebras generating $A$. By Lemma 4.1, it follows  that for large enough 
$n$ we have $\|\psi_k\circ E_{A_n'\cap M}\|< \varepsilon/2$, $\forall k$. Thus, if we let $\{q_i\}_i$ denote the minimal projections of $A_n$ 
then for any $\varphi_k$ in $F_0$ we have:
$$
\|\sum_i q_i \cdot \varphi_k \cdot q_i\| \leq \|\varphi_k - \psi_k\| + \|\sum_i q_i \cdot \psi_k \cdot q_i\|< \varepsilon  \tag 5.4.1
$$ 

If $M$ is II$_1$, one can obviously slightly perturb the partition $\{q_i\}_i$ so that each $q_i$ has trace 
of the form $k_i/2^t$, for some integers $t\geq 1$, $k_i\geq 1$, while still having the (strict!) inequality $(5.4.1)$ satisfied. 
Then any refinement of this partition  to a partition of 1 with $2^t$ projections of trace $2^{-t}$, 
$\{p_j\}_j \subset \Cal P(M)$, will do. 

If $M$ is properly infinite, then by Lemma 4.6 we can slightly perturb $\{q_i\}_i$ in the $s$-topology to projections that are all infinite (thus all $\sim 1$) 
and still satisfy the inequality $(5.4.1)$. 
\hfill $\square$ 

\vskip.05in

We can now prove the main result of this section, showing that any MASA in any continuous separable factor is flat. We in fact also clarify the 
remaining atomic  (type I) case, where the opposite phenomenon occurs, as  we'll show that no MASA can be flat.

\proclaim{5.5. Theorem} $1^\circ$ If $M$ is a separable factor of type $\text{\rm II}$ or $\text{\rm III}$ $($i.e., if $M$ is continuous$)$, then any MASA in $M$ is flat. 

$2^\circ$ If $M$ is a separable, non-trivial factor of type $\text{\rm I}$ $($i.e., if $M$ is atomic$)$, then all MASAs in $M$ are non-flat.  
\endproclaim
\noindent
{\it Proof of} $1^\circ$ Assume first that $M$ is of  type II$_1$, with its (unique) normalized trace $\tau$. 
Let $A\subset M$ be a MASA, $\Cal F \subset L^1M$ a finite set of positive elements 
of trace $1$ and $\delta >0$. We have to prove that there exists $u\in \Cal U(M)$ such that $\text{\rm sup}\{|\tau(uau^*(\xi - \zeta))| \mid  a\in (A)_1\} 
< \delta$, $\forall \xi, \zeta \in \Cal F$. Without loss of generality, we may clearly assume $\Cal F\subset M$ (i.e., 
that all elements in $\Cal F$ are bounded). 
But this is (Lemma 2.3 in [P17]), applied to $P=M$, $Q=A$ and $F=\{\xi - \zeta \mid \xi, \zeta \in \Cal F\}$. 

Thus, if $M$ is II$_1$ then any MASA in $M$ is flat. 

Assume now that $M$ is of type II$_\infty$ and let $A\subset M$ be a MASA that's generated by projections of finite trace $Tr=Tr_M$. 
To show that $A$ is flat in $M$, we have to prove that for any finite set $\Cal F \subset L^1(M, Tr)$ of positive elements of trace $Tr$ equal to $1$ and any $\varepsilon >0$, 
there exists $u\in \Cal U(M)$ such that $\text{\rm sup}\{|Tr(uau^*(\xi - \zeta))| a\in (A)_1\} 
< \varepsilon$, $\forall \xi, \zeta \in \Cal F$. By small $\| \ \|_{1,Tr}$-perturbation of all elements in $\Cal F$, it is clearly sufficient to 
show this for finite sets $\Cal F$ with all its elements bounded in operator norm that are supported by projections of finite  trace $Tr$. This latter assumption implies 
that we may suppose there exists a finite projection $p\in M$ such that $p\Cal F p=\Cal F$. 

Since $A$ is generated by finite projections and it is diffuse, there exists $e\in \Cal P(A)$ such that $Tr(e)=Tr(p)$. Let $u_0\in \Cal U(M)$ 
be so that $u_0eu_0^*=p$. By replacing $A$ with $u_0Au_0^*$, we may assume the support $p$ of $\Cal F$ lies in $A$, $p\in A$. 
We can apply now the first part of the proof 
to the II$_1$ factor $pMp$ with its MASA  $Ap\subset pMp$ and to $\delta=\varepsilon/Tr(p)$, to get a unitary element $v\in pMp$ 
such that $\text{\rm sup}\{|\tau_{pMp}(vav^*(\xi - \zeta))| a\in (Ap)_1\} 
< \delta$, $\forall \xi, \zeta \in \Cal F$. But since $\xi, \zeta$ are supported by $p\in A$ and $\tau_{pMp}(x)=Tr(x)/Tr(p)$ 
for any $x\in pMp$, if we denote $u=v+(1-p)$ then we have 
$$
\text{\rm sup}\{|Tr(uau^*(\xi - \zeta))| a\in (A)_1\}/Tr(p)
$$
$$
=\text{\rm sup}\{|\tau_{pMp}(vav^*(\xi - \zeta))| a\in (Ap)_1\}  < \varepsilon/Tr(p)
$$
for all $\xi, \zeta \in \Cal F$. 

This ends the proof of the fact that $A$ is flat in $M$. Combining  with the first part and using 
Lemma 5.4, it follows that any MASA in any factor of type II is flat. 

The case when $M$ is of type III$_\lambda$, $0<\lambda <1$, follows now 
by taking into account that such a factor contains a MASA $A\subset M$ with the property that there exists an intermediate subfactor $A\subset N \subset M$ 
of type II$_1$ and then using the II$_1$ case above in combination with Lemma 5.3. 

To settle the remaining cases $M$ of type III$_1$ or type III$_0$, let us first show that any AFD type III$_0$ or III$_1$ factor contains flat MASAs. 

Indeed, if $M$ is the unique (by [H87], [C85]) AFD factor of type III$_1$, $M=\Cal R_1$, 
then by [H87] it contains an ergodic copy of $R$ with a normal conditional 
expectation. By [P81a], this implies $R$ contains a diagonal $D\subset R$ that's a MASA in $\Cal R_1$. Thus, since $D$ is flat in $R$, 
by Lemma 5.3 it is flat in $\Cal R_1$. 

Take now $M$ to be an AFD factor of type III$_0$, $M = \Cal R_0$. Let $D\subset \Cal R_0$ be its (unique up to conjugation by automorphisms) Cartan subalgebra 
and $\phi$ a normal faithful state on $\Cal R_0$ that has $D$ in its centralizer. To prove flatness of MASAs in $\Cal R_0$ it is sufficient to 
check condition $3^\circ$ of Lemma 5.4. Also, it is sufficient to consider finite sets $F_0$ 
in a total subset $\Cal L$ of $\{\varphi \in L^1(\Cal R_0, \phi) \mid \varphi(1)=0\}$. We'll take $\Cal L$ to be union between the set $\Cal L_1$  of elements of the form 
$\phi( \cdot \ ua)$ with $a\in D$ and $u$ in the part of the normalizing groupoid of $D$ that's outer on $D$, union with the set $\Cal L_0=\{\phi( \cdot \ b) \mid b \in D, \phi(b)=0\}$. 
By Lemma 5.4, we need to construct a partition of $1$ with projections in $M$ that ``kill'' the union between two given finite sets of elements $F_0\subset \Cal L_0$ and  $F_1 \subset \Cal L_1$. 

One can obviously first find a partition $\{p_i\}_i$ of $1$ in $D$ that kills $F_1$ and such that each $p_i$ is $\phi$-independent to $F_0$ (exercise!). 
By approximating $bp_i \in p_iF_0p_i$ by elements in a finite dimensional subalgebra of $Dp_i$ 
with minimal projections equivalent in $M$, in each $p_iMp_i$ we can find partitions of $1$ that ``kill'' all $F_0p_i$ (e.g., like in [P81c]). This shows 
that the AFD III$_0$ factor $\Cal R_0$ has flat MASAs. 

Thus, at this point we know that any MASAs in any continuous AFD factor is flat. 
 
Assume now that $M$ is an arbitrary type III$_1$ factor. By [P84], $M$ contains an irreducible AFD subfactor $\Cal R \subset M$ with the property that $\Cal R$ contains a 
MASA $A$ of $M$. Since $A$ is flat in $\Cal R$, it is also flat in $M$ by Lemma  5.3. Thus, all MASAs in $M$ are flat. 

If in turn $M$ is III$_0$, then by [P83] it contains an AFD III$_0$ subfactor $\Cal R_0 \subset M$ with a Cartan subalgebra $D\subset \Cal R_0$ that's a MASA 
in $M$. But $D$ is flat in $\Cal R_0$ so by lemma 5.3 it is also flat in $M$. 

\vskip.05in 
\noindent
{\it Proof of} $2^\circ$ If $M=\Cal B(\ell^2_n)$ for some finite $n\geq 2$, then let $F$ be a finite set of states S$_{\text{\rm n}}(M)$ on $M$ 
with the property that any state on $M$ is at distance $\leq 1/2$ from a state in $F$ (i.e., $F$ is $1/2$-dense in S$_{\text{\rm n}}(M)$). 
If $A\subset M$ is any arbitrary MASA in $M$, then 
the restrictions of $F$ to $A$ give a $1/2$-dense set of states on $A$. Take now $\varphi, \psi \in \text{\rm S}_{\text{\rm n}}(M)$ with disjoint supports in $A$ and let   
$\varphi', \psi'\in F$ be so that $\|(\varphi-\varphi')_{|A}\|<1/2$, $\|(\psi-\psi')_{|A}\|<1/2$. Then $\|\varphi'-\psi'\|\geq \|\varphi - \psi\|-\|\varphi - \varphi'\| 
-\|\psi - \psi'\|\geq 1$. Thus, given any unitary conjugate $A=uDu^*$ of the diagonal $D\subset M$, there exist two states in $F$ that are ``far apart'' when restricted to $A$, 
showing that $D$  (and any conjugate of it) is not flat. 

Let now $M=\Cal B(\ell^2\Bbb N)$. By Lemma 5.3 it is sufficient to show that its diffuse MASA is non-flat.  We view this MASA as $A=L^\infty([0,1], \mu) \subset \Cal B(L^2([0,1], \mu))=M$, 
where $\mu$ is the Lebesgue measure. 
Let $\xi_0=2^{1/2}\chi_{[0,1/2]}$, $\xi_1=2^{1/2}\chi_{[1/2,1]}$, $\xi=1=\chi_{[0,1]}$. If we assume $A$ is flat, then one can choose 
$U\in \Cal U(L^2([0,1]))$ such that the unit vectors 
$\eta_0=U(\xi_0)$, $\eta_1=U(\xi_1)$, $\eta=U(\xi)$ implement ``almost the same'' measure on $[0,1]$. 
In other words, $|\eta_0|^2, |\eta_1|^2, |\eta|^2$  are almost equal as elements in $L^1([0,1])$. Equivalently, $|\eta_0|, |\eta_1|, |\eta |$ 
are almost equal in $L^2([0,1])$. By multiplying  $U$ on the left by a unitary in $L^\infty([0,1])$ that agrees with the phase of $\eta_0$, 
we may assume $\eta_0\geq 0$.   

Altogether, the unit vectors $\eta_0, \eta_1, \eta\in L^2([0,1])$ have almost equal absolute value, with $\eta_0\geq 0$, $\eta_0 \perp \eta_1$, 
and $\eta=2^{-1/2}\eta_0 + 2^{-1/2}\eta_1$. Since $\eta\eta^*=\eta_0^2/2+ \eta_1\eta_1^*/2 + \eta_0\eta_1^*$, this implies $\|\eta_0  \eta_1^*\|_1\approx 0$. 
Equivalently, $\int |\eta_0| |\eta_1| \text{\rm d}\mu\approx 0$. 
Thus, $0\approx \| |\eta_0|^2 - |\eta_1|^2\|_1 \geq \||\eta_0|-|\eta_1| \|_2^2 \approx  2$, a contradiction. 

\hfill $\square$

\vskip.05in 
\noindent
{\bf 5.6. Remark}. Note that the fact that the diffuse (and thus all) MASAs in $\Cal B(\ell^2\Bbb N)$ are not flat follows also as a consequence 
of Theorem 6.4. More precisely, of its poof, which shows that if  $M\subset \Cal M$ is an inclusion of separable factors so that $M$ contains a diffuse MASA of $\Cal M$ 
that's flat  in $M$, then $M$ contains a copy of the hyperfinite II$_1$ factor $R\subset M$ that has trivial relative commutant in $\Cal M$. Thus, it 
does apply to the case $M=\Cal M=\Cal B(\ell^2\Bbb N)$, if we assume the diffuse MASA in $\Cal B(\ell^2\Bbb N)$ is flat. But by von Neumann's Bicommutant Theorem ([vN29]), 
$R'\cap \Cal B(\ell^2\Bbb N)=\Bbb C$ implies $\Cal B(\ell^2\Bbb N)=R''=R$, 
contradicting the fact that $R$ is type II while $\Cal B(\ell^2\Bbb N)$ is type I ([MvN36]).

\heading 6. From MASA-ergodicity to $R$-ergodicity
\endheading 

In this section we prove the main result in this paper, showing the implication ``MASA-ergodicity $\Rightarrow$ $R$-ergodicity'' for an inclusion of separable factors. 

We begin with a criterion for a representation of a UHF algebra to give rise to the hyperfinite II$_1$ factor. This can be derived from results 
(Section 2 in [P67]), but we have included a self-contained argument, for the readers's convenience. 

\proclaim{6.1. Lemma}  Let $B_0$ be a  UHF algebra, obtained as the $C^*$-inductive limit of matrix factors $B_{0,n}\simeq \Bbb M_{k_n}(\Bbb C)$ 
with $k_n | k_{n+1}$, $\forall n$, and denote by $\tau$ its unique trace state. 
Let $B_0\subset \Cal B(\Cal H)$ be a representation of $B_0$ on a separable Hilbert space $\Cal H$ and $\{\xi_n\}_n\subset \Cal H$ 
a sequence of unit vectors in $\Cal H$ that's dense in the set of unit vectors of $\Cal H$. The following conditions are equivalent 

\vskip.05in
$1^\circ$ The von Neumann algebra $B_0''=\overline{B_0}^{wo}$ is the hyperfinite $\text{\rm II}_1$ factor $R$; 

\vskip.03in
$2^\circ$ Given any unit vector $\xi \in \Cal H$, the state $\varphi$ on $B_0$ implemented by $\xi$ satisfies 
$\underset{n \rightarrow \infty}\to{\lim} \|\tau_{|B_{0,n}'\cap B_0}-\varphi_{|B_{0,n}'\cap B_0}\|=0$. 

\vskip.03in
$3^\circ$ For any $m$, the state $\varphi_m$ on $B_0$ implemented by $\xi_m$ satisfies 
$\underset{n \rightarrow \infty}\to{\lim} \|\tau_{|B_{0,n}'\cap B_0}-{\varphi_m}_{|B_{0,n}'\cap B_0}\|$ $=0$.

\endproclaim
\noindent
{\it Proof}. $1^\circ \Rightarrow 2^\circ$ If $B_0''=R$ is the hyperfinite II$_1$ factor, then $\Cal H$ is a direct sum of cyclic Hilbert $R$-modules 
of the form $L^2Rp$, with $p\in \Cal P(M)$. So it is clearly sufficient to prove $2^\circ$ for $\Cal H$ itself of this cyclic form. 
But then any vector state $\varphi$ on $R\subset \Cal B(L^2Rp)$ is of the form $\varphi(x)=\tau(xb)$ for some $b \in L^1R_+$. 

If $E_{B_{0,n}}$ denotes the unique expectation of $R$ onto $B_{0,n}$ that preserves the tracial state $\tau$ on $R$, then $\lim_n \|E_{B_{0,n}}(b)-b\|_1=0$. 
Thus, for any $\varepsilon >0$ there exists $n_0$ such that if $n\geq n_0$ then $\|E_{B_{0,n}}(b)-b\|_1< \varepsilon.$
Since $B_{0,n}$ are factors, if $b_0\in B_{0,n}$ then $E_{B_{0,n}'\cap R}(b_0)=\tau(b_0)1$. Thus, 
if we take $b_0=E_{B_{0,n}}(b)$, then for any $x\in (B_{0,n}'\cap R)_1$ and any $n\geq n_0$  we have 
$$
|\tau(x)-\varphi(x)|=|\tau(x)-\tau(xb)| < |\tau(x)-\tau(xE_{B_{0,n}}(b))|+\varepsilon 
$$
$$
=|\tau(x)-\tau(xE_{B_{0,n}'\cap R}(b_0))|+\varepsilon = |\tau(x)-\tau(b_0)\tau(x)|+\varepsilon = \varepsilon. 
$$

$2^\circ\Rightarrow 3^\circ$ is trivial. To prove $3^\circ \Rightarrow 1^\circ$, denote $\Cal R=B_0''\subset \Cal B(\Cal H)$. Since the UHF algebra $B_0$ 
has Dixmier's property and unique trace, the finite part of $\Cal R$ is isomorphic to $R$. So if  
$\Cal R\not\simeq R$, then $\Cal R$ has a non-zero properly infinite part. Let $\xi\in \Cal H$ be a unit vector with $p'=[B_0\xi] = [\Cal R\xi] \in \Cal R'$ 
having the property that $\Cal Rp'$ is properly infinite. 

Let $m$ be so that $\|\xi-\xi_m\| < 2^{-5}$ and $n$ be so that 
$|\tau_{|B_{0,n}'\cap B_0}-{\varphi_m}_{|B_{0,n}'\cap B_0}\|<2^{-5}$.  Thus, $\xi'=p'(\xi_m)$ satisfies $\|\xi'-\xi_m\|<2^{-4}$, $\|\xi'\|\geq 15/16$, while the vector functional $\varphi'$ 
implemented by $\xi'$ satisfies $\|\varphi'_{|B_{0,n}'\cap B}-{\varphi_m}_{|B_{0,n}'\cap B_0}\|<2^{-3}$. Taken together, these conditions imply 
$\||\tau_{|B_{0,n}'\cap B_0}- {\varphi'}_{|B_{0,n}'\cap B_0}\|<2^{-2}$. Thus, the vector  $\xi'\in \Cal H'=\overline{\Cal R\xi}$, which has norm $2^{-4}$-close to $1$, implements on the properly 
infinite von Neumann algebra $\Cal R_0=(B_{0,n}'\cap \Cal R)p'\subset \Cal B(\Cal H')$ a positive normal functional with the property that on the $s^*$-dense 
subalgebra $B_{0,n}'\cap B$ is $2^{-2}$-close to the trace state $\tau$. It follows that for any $u\in \Cal U(B_0)$ and $x\in (\Cal R_0)_1$ we have 
$|\varphi'(uxu^*)-\varphi'(x)| \leq 2^{-2}$. 

Taking an amenable subgroup $\Cal U_0\subset \Cal U(B_0)$ with the property that the span of $\Cal U_0$ is norm-dense in $B_0$ 
and ``integrating'' over $\Cal U_0$ with respect to an invariant mean on $\Cal U_0$ in the unit ball of $\Cal R_0^*$ (which is $\sigma(\Cal R_0^*, \Cal R_0)$-compact), 
we get a positive functional $\varphi_0$ on $\Cal R_0$ that's $\Cal U_0$-invariant and is still $2^{-2}$-close to $\varphi'$, which in turn is positive normal with $\varphi'(1)\geq (\frac{15}{16})^2 \geq 7/8$. 
So if we denote by $\psi_0$ the normal part of $\varphi_0$, then $\psi_0$ is still $\Cal U_0$-invariant and $2^{-2}$-close to $\varphi'$ (and therefore $\psi_0\neq 0$). 
But $\Cal U_0''=\Cal R_0$, so this implies $\psi_0$ is a non-zero positive normal tracial functional on $\Cal R_0$, contradicting the fact that $\Cal R_0$ is properly infinite. 
\hfill 
$\square$

\proclaim{6.2. Lemma} Let $M \subset \Cal M$ be a MASA-ergodic embedding of continuous factors. Given any finite dimensional factor 
$Q\subset M$, any finite set $F_0$ of normal states on $\Cal M$ and any $\delta_0 >0$, there exists an abelian von Neumann 
subalgebra $A\subset Q'\cap M$ which is a MASA in $Q'\cap \Cal M$ and satisfies the conditions:

$$\|(\varphi-\psi)_{|A}\|< \delta_0, \forall \varphi, \psi \in F_0. \tag 6.2.1  $$

$$\|\varphi_{|Q \vee A}-\varphi_{|Q}\otimes \varphi_{|A}\|< \delta_0, \forall \varphi \in F_0. \tag 6.2.2  $$

\endproclaim
\noindent
{\it Proof}. Let $\{e_{ij}\}_{1\leq i,j \leq n}$ be a set of matrix units for $Q$. 
Since $M\subset \Cal M$ is MASA-ergodic, there exists an abelian von Neumann subalgebra $B\subset M$ that's a MASA in $\Cal M$. After some  
unitary conjugation, we may clearly assume $e_{11}\in Q$.  Thus, $Be_{11}\subset e_{11}Me_{11}$  is a MASA in $e_{11}\Cal M e_{11}$. 
This implies $A_0=\{\sum_i e_{i1}be_{1i} \mid b\in e_{11}B\} \subset Q'\cap M$ is a MASA in $Q'\cap \Cal M$, which is flat by Theorem 5.5 

By applying flatness to the finite set of normal functionals $F=\{ (\varphi - \psi)_{|Q'\cap M} \mid \varphi, \psi \in F_0\} \cup \{ ((\varphi (e_{ij} \ \cdot) - \varphi(e_{ij})\varphi)_{|Q'\cap M} \mid 
\varphi \in F_0, 1\leq i,j \leq n \}$, which all vanish at $1$, it follows that for  $\delta=\delta_0/n^2$, there exists a unitary element in $Q'\cap M$ such that $A=uA_0u^*$  
and $\varphi, \psi \in F_0$ satisfy

$$\|(\varphi-\psi)_{|A}\|< \delta_0,  $$

$$|\varphi (e_{ij}a) -\varphi(e_{ij})\varphi(a)| < \delta_0/n^2, 1\leq i, j \leq n, a\in (A)_1 $$

But the second condition  clearly implies $(6.2.2)$ while the first condition is just $(6.2.1)$. 
\hfill $\square$

\proclaim{6.3. Lemma} Let $M$ be an infinite dimensional factor, $F$ a finite family of normal functionals on $M$   
and $\{e_{i}\}_{1\leq i \leq n}$ a partition of $1$ with mutually equivalent projections in $M$. Given any $\delta >0$ there exists a set of matrix 
units $\{e_{ij}\mid 1\leq i, j \leq n\} \subset M$ such that  $e_{ii}=e_i$ and $|\varphi(e_{ij})|\leq \delta$, $\forall i\ne j$, $\forall \varphi \in F$. 
\endproclaim
\noindent
{\it Proof}. We construct by induction over $k=1, 2, ..., n$ some matrix units $\{e_{ij}\mid 1\leq i,j \leq k\}\subset M$,  
such that $e_{ii}=e_i$, $1\leq i \leq k$, and $|\varphi(e_{ij})|< \delta$, $\forall 1\leq i\neq j \leq k$.  

Assume we have done this up to some $k<n$. Choose a partial isometry $v\in M$ such that $v^*v=e_{k+1}$ and $vv^*=e_1$ 
and consider the finite set of normal functionals on $e_{k+1}Me_{k+1}$ given by 

$$
F_1=\{\varphi(e_{j1}v  \ \cdot \ ) \mid \varphi \in F, 1\leq j \leq k \} \cup \{ \varphi ( \cdot \ v^* e_{1j}) \mid \varphi \in F, 1\leq j \leq k\}
$$ 

Take a copy of the diffuse abelian von Neumann algebra $L^\infty(\Bbb T)$ on $e_{k+1}Me_{k+1}$, which we view as generated 
by a Haar unitary $u$. Since $u^m$ tends weakly to $0$ as 
$|m|\rightarrow \infty$, there exists $m$ large enough such that $|\psi(u^{\pm m})| < \delta$ for all  $ \psi \in F_1$. Thus, if for each 
$1\leq j \leq k$ we define $e_{j,k+1}=e_{j1}vu^m$ and  
$e_{k+1,j}=u^{-m}v^*e_{1j}$, then all conditions are satisfied. 

\hfill $\square$ 

\proclaim{6.4.  Theorem} Let $M\subset \Cal M$ be an embedding of continuous, separable factors. 
If $M \subset \Cal M$ is MASA-ergodic, then it is $R$-ergodic. Moreover, there exists an embedding $R\hookrightarrow M$ 
such that the diagonal $D$ of $R$ is a MASA in $\Cal M$. If in addition $\Cal M\subset \Cal B(\Cal H)$ is a normal representation,   
$\xi\in \Cal H$ is a separating unit vector for $\Cal M$, and $\varepsilon >0$, then $R\subset M$ can be chosen so 
that the state $\phi$  implemented by $\xi$ on $R$ is equal to the trace $\tau$ of $R$ when  restricted to $D$, and satisfies $\|\phi_{|R}-\tau\|<\varepsilon$. 
\endproclaim
\noindent
{\it Proof}. Since the case $\Cal M$ is of type II$_1$ of this theorem was proved in [P81a], we may and will assume from now on that $\Cal M$ is a 
properly infinite factor. 

The reader should notice that all we will use for the rest of the proof is the assumption that $M$ contains a diffuse flat MASA of $\Cal M$. 
Thus, the arguments that follow do also apply to the case $M=\Cal M=\Cal B(\ell^2\Bbb N)$, when assuming that the diffuse MASA 
in $\Cal B(\ell^2\Bbb N)$ is flat (see Remark 5.6).

Let $\Cal M\subset \Cal B(\Cal H)$ be the standard representation of $\Cal M$ on a separable Hilbert space. 
Let $\xi\in \Cal H$ be a cyclic and separating unit vector for $\Cal M$, which in the case $M$ is II$_1$ 
we assume to be the trace when restricted  to $M$. To see that this is possible, notice that if $\eta\in \Cal H$ is a cyclic and separating unit 
vector for $\Cal M$ and $b\in L^1M_+$ is so that $\langle x \eta, \eta\rangle =\tau(xb)$, $\forall x\in M$, 
and we let $p_n$ be the spectral projection of $b$ corresponding to  the interval $[n, \infty)$, 
then $\{p_nb^{-1/2}\eta\}_n$ is Cauchy in $\Cal H$ and its limit $\xi \in \Cal H$ is a unit vector 
which is still cyclic and separating for $\Cal M$ while $\langle x\xi, \xi\rangle=\tau(x)$, $\forall x\in M$. 

We denote $\phi=\omega_\xi$ the normal faithful state implemented by $\xi$ on $\Cal M$. 
Let also $\{\xi_n\}_{n\geq 1} \subset  \Cal H$  be a sequence of unit vectors that's dense in the set of unit vectors in $\Cal H$ 
and with $\xi_1=\xi$. 
Denote by $\varphi_j$ the state implemented by $\xi_j$ on $\Cal M$ (so $\varphi_1=\phi$). 

We construct recursively an increasing sequence of dyadic factors $Q_m$ in $M$, 
with matrix units $\{e^m_{ij} \mid 1\leq i,j \leq 2^{k_m}\}$, with $Q_0=\Bbb C1$, such  that 
if we denote $D_m=\sum_i \Bbb C e^m_{ii}$, then for each $m$ there exist elements $\{x^m_i\}_i \subset (\Cal M')_1$ satisfying the properties: 

$$
\| [D_m \xi](\xi_j) - \sum_i e^m_{ii} x^m_i (\xi_j)\| <  2^{-m},  1\leq j \leq  m. \tag 6.4.1
$$

$$
\|(\varphi_k - \tau)_{|Q_{k-1}'\cap Q_m}\| < 2^{-k}, 1\leq k \leq m.  \tag 6.4.2 
$$

Assume we have constructed these algebras up to $m=n$. We then apply Lemma 6.2 to $Q=Q_n$ and an arbitrarily 
small $\delta>0$ to get an abelian subalgebra $A\subset Q_n'\cap  M$ that's a MASA in $Q_n'\cap \Cal M$ and satisfies

$$
\|(\varphi_i-\varphi_j)_{|A}\|< \delta/4,  1\leq i,j \leq n+1. \tag 6.4.3  
$$

$$
\|{\varphi_j}_{|Q_n \vee A}-{\varphi_j}_{|Q_n}\otimes {\varphi_j}_{|A}\|< \delta, 1\leq j \leq n+1. \tag 6.4.4  
$$

Let $A_p\subset A$, $p\geq 1$, be an increasing sequence of finite dimensional von Neumann subalgebras that exhaust $A$. Thus, if we denote $D_n=\text{\rm sp}\{e^n_{jj}\}_j$, then 
$B=A\vee D_n$ is a MASA in $\Cal M$ and $B_p=A_p \vee D_n \nearrow B$. By Lemma 4.3, it follows that there exists $p_0$ such that if we denote $\{f_k\}_k$ the minimal 
projections of $B_{p_0}=D_n \vee A_{p_0}$ then there exist $\{x_k\}_k \subset (\Cal M')_1$ with the property that 
$$
\|[B_{p_0}\xi](\xi_i)-[B\xi](\xi_i)\| + \|[B\xi](\xi_i)-\sum_k f_k x_k(\xi_i)\| < \delta, 1\leq i \leq n+1. 
$$

By Lemma 4.4,  if $A_{p_0}\subset A^0\subset A$ is an arbitrary  intermediate finite dimensional von Neumann subalgebra and we denote  
by $\{f_j^0\}_j$ the minimal projections of $B^0:=D_n \vee A^0$ and by $\{x^0_j\}_j$ an appropriate rearrangement with repetition of $\{x_k\}_k$, then we still have 
$$
\|[B^0\xi](\xi_i)-[B\xi](\xi_i)\| + \|[B\xi](\xi_i)- \sum_j f_j^0 x_j^0(\xi_i)\| < \delta, 1\leq i \leq n+1. \tag 6.4.5
$$

Note at this point that if $M$ is a II$_1$ factor, then the fact that $\Cal M$ is properly infinite automatically implies that all $f_j^0$ are 
infinite in $\Cal M$. On the other hand, if $M$ is a properly infinite factor, then by Lemma 4.5 we may slightly perturb the finite partition given by $A^0$ to assume all of its 
minimal projections are infinite in $Q_n'\cap \Cal M$ (so in $\Cal M$ as well) and have $\phi$-value of the form $\frac{k_j}{2^p}$, on all of 
its minimal projections, while the resulting $B^0, \{f_j^0\}_j$ and $A^0$ still satisfy $(6.4.5)$. 

By Lemma 4.6, we can then refine $A^0$ to a finite dimensional partition $A^1$ contained in $A$ that has all its minimal projections infinite (thus equivalent) in $\Cal M$ and of $\phi$-value $\frac{1}{2^p}$. 
If we denote $D_{n+1}$ the finite dimensional algebra $D_n \vee A^1$ and $\{e^{n+1}_{jj}\}_j$ its minimal projections, then by Lemma 4.4 
we have 
 $$
\|[D_{n+1}\xi](\xi_l)-B\xi](\xi_l)\| + \|[B\xi](\xi_l)- \sum_j e^{n+1}_{jj}x^{n+1}_j(\xi_l)\| < \delta, 1\leq l \leq n+1, 
$$
where $\{x^{n+1}_j\}_j\subset (\Cal M')_1$ is an appropriate rearrangement with repetition of the elements $\{x_j^0\}_j$ appearing in $(6.4.5)$. 
By the triangle inequality,  this implies 
 $$
\|[D_{n+1}\xi](\xi_l)- \sum_j e^{n+1}_{jj}x^{n+1}_j(\xi_l)\| < \delta, 1\leq l \leq n+1, \tag 6.4.6
$$
which if we take $\delta\leq 2^{-n-1}$ amounts to  the inequality $(6.4.1)$ for $m=n+1$.

By Lemma 4.6, we can now take a system of matrix units $\{e^{n+1}_{ij}\}_{i,j}\subset M$, having $e^{n+1}_{jj}$ as the diagonal, so that 
for all $1\leq i\neq j\leq 2^{k_{n+1}}$ we have 
$$
|\varphi_l(xe^{n+1}_{ij})| < 2^{-2k_{n+1}} \delta/2, \forall x\in (Q_n)_1, 1\leq l \leq n+1  \tag 6.4.7 
$$

We claim that if we take $\delta< 2^{-n-1}$ in $(6.4.7)$, then the inequalities $(6.4.2)$ are satisfied as well. To see this, 
let $y$ be an arbitrary element in the unit ball of $Q_{k-1}'\cap Q_{n+1}=(Q_{k-1}'\cap Q_n)\vee (Q_n'\cap Q_{n+1})$, which we write it in the form 
$y=\sum_{i,j} y_{ij} e^{n+1}_{ij}$ with $y_{ij}\in (Q_{k-1}'\cap Q_n)_1$. Taking into account that for  $i\neq j$ we have $\tau(e^{n+1}_{ij})=0$ and 
$|\varphi_k(y_{ij}e^{n+1}_{ij})|\leq 2^{-2k_{n+1}}\delta$ (due to $(6.4.7)$), it follows that: 
$$
|(\varphi_k-\tau)(y)| = |\sum_{i,j} (\varphi_k-\tau)(y_{ij}e^{n+1}_{ij})|  \tag 6.4.8
$$
$$
\leq \sum_{i\neq j} |(\varphi_k-\tau)(y_{ij}e^{n+1}_{ij})| + |\sum_{j} (\varphi_k-\tau)(y_{jj}e^{n+1}_{jj})|
$$
$$
=\sum_{i\neq j} |\varphi_k(y_{ij}e^{n+1}_{ij})| + |\sum_{j} (\varphi_k-\tau)(y_{jj}e^{n+1}_{jj})|
$$
$$
\leq 2^{2k_{n+1}}2^{-2k_{n+1}} \delta/2 + |\sum_{j} (\varphi_k-\tau)(y_{jj}e^{n+1}_{jj})| 
$$
$$
= \delta/2  + |\sum_{j} (\varphi_k-\tau)(y_{jj}e^{n+1}_{jj})|. 
$$
For this last term, we get the estimate
$$
\delta/2+ |\sum_{j} (\varphi_k-\tau)(y_{jj}e^{n+1}_{jj})| = \delta/2+ |(\varphi_k-\tau)(\sum_{j} y_{jj}e^{n+1}_{jj})| 
$$
$$
\leq \delta/2+ |\varphi_k(\sum_j y_{jj}e^{n+1}_{jj})-\sum_j \varphi_k(y_{jj})\varphi(e^{n+1}_{jj})| 
$$
$$
+  \sum_j |\varphi_k(y_{jj})\varphi_k(e^{n+1}_{jj})-\tau(y_{jj})\tau(e^{n+1}_{jj})| 
$$
which by first using $(6.4.4)$ and then the triangle inequality  we majorize by 
$$
\delta/2 + \delta/4+\sum_j |\varphi_k(y_{jj})\varphi_k(e^{n+1}_{jj})-\tau(y_{jj})\tau(e^{n+1}_{jj})| 
$$
$$
\leq 3\delta/4 + \sum_j |\varphi_k(y_{jj})| |\varphi_k(e^{n+1}_{jj})-\tau(e^{n+1}_{jj})| + \sum_j |\varphi_k(y_{jj})-\tau(y_{jj})| \tau(e^{n+1}_{jj}).  
$$
Since by $(6.4.3)$ we have $|\varphi_k(y_{jj})-\tau(y_{jj})| \leq \delta/4$, while the fact that $\tau(e^{n+1}_{jj})=2^{-k_{n+1}}$ and  $|\varphi_k(y_{jj})|\leq 1$, $\forall j$, 
entails $\sum_j |\varphi_k(y_{jj})-\tau(y_{jj})| \tau(e^{n+1}_{jj}) \leq \delta/4$, it follows that the above last term is majorized by 
$$
3\delta/4+ \sum_j|\varphi_k(e^{n+1}_{jj})-\tau(e^{n+1}_{jj})|.  
$$

If we now denote by $b_k \in L^1A_+$ the Radon-Nykodim derivative of $\varphi_k$ with respect to $\tau=\varphi_1$, then 
by $(6.4.3)$ we have $\|b_k-1\|_1 \leq \delta/4$, so the above last term can be further majorized by  
$$
3\delta/4+ \sum_j|\varphi_k(e^{n+1}_{jj})-\tau(e^{n+1}_{jj})| = 3\delta/4 + \sum_j |\varphi_k((b_k-1)e^{n+1}_{jj})| 
$$
$$
\leq 3\delta/4 + \sum_j \|(b_k-1)e^{n+1}_{jj}\|_1 =3\delta/4 + \|b_k-1\|_1 \leq \delta,  
$$
thus finalizing our estimate $(6.4.8)$. Since $y\in  (Q_{k-1}'\cap Q_{n+1})_1$ was arbitrary in all this, it follows that 
$$
\|(\varphi_k-\tau)_{|Q_{k-1}'\cap Q_{n+1}}\| = \sup\{|(\varphi_k-\tau)(y)| \mid y\in (Q_{k-1}'\cap Q_{n+1})_1\} \leq \delta, \tag 6.4.9
$$
which if we take $\delta< 2^{-n-1}$ shows that $Q_{n+1}$ satisfies  $(6.4.2)$ for $m=n+1$.

With the $Q_m, D_m$ constructed recursively for $m\geq 1$ so that to satisfy $(6.4.1)$, $(6.4.2)$, let us now define 
$D=\overline{\cup_m D_m}^w$. By condition $(6.4.1)$ and Lemma 4.3, it follows that $D$ is a MASA in  
$\Cal M$. Also, if we define $R_0=\overline{\cup_m Q_m}^n$, then by Lemma 6.1 and condition $(6.4.2)$ 
it follows that the weak operator closure $R$ of $R_0\subset \Cal B(\Cal H)$ 
is the hyperfinite II$_1$ factor. Since $R$ contains $D$, it also follows that $R'\cap \Cal M = R'\cap (D'\cap \Cal M)=R'\cap D=\Bbb C$, 
showing that $M\subset \Cal M$ is $R$-ergodic. 

\hfill $\square$

\proclaim{6.5. Corollary} Any separable factor $\Cal M$  that's not of type I admits an ergodic embedding of $R$. 
\endproclaim

\noindent
{\it Proof}. This is just the case $M=\Cal M$ of Theorem 6.4. 
\hfill $\square$

\proclaim{6.6. Corollary} Any continuous separable factor $\Cal M$  
can be embedded ergodically into the unique separable AFD type $\text{\rm II}_\infty$ factor $R^\infty=R \overline{\otimes} \Cal B(\ell^2\Bbb N)$. 
\endproclaim

\noindent
{\it Proof}. By Corollary 6.5, there exists an ergodic embedding $R\hookrightarrow \Cal M^{op}$. Thus, if $\Cal M \subset \Cal B(L^2\Cal M)$ 
is the standard representation of $\Cal M$, with the corresponding embedding $R\subset \Cal M'\simeq \Cal M^{op}$, 
then one gets an embedding $\Cal M \subset R'\cap \Cal B(L^2\Cal M) \simeq R^\infty$, which satisfies 
$\Cal M'\cap R^\infty=\Cal M'\cap (R'\cap \Cal B(L^2\Cal M))=R'\cap \Cal M^{op}=\Bbb C$. 
\hfill $\square$

\vskip.05in

As we mentioned before, a class of ergodic embeddings $M \subset \Cal M$ that are particularly interesting are the 
II$_\infty$ $\subset$ III$_1$ inclusions of factors arising from the continuous decomposition of a type III$_1$ factors $\Cal M$. By [H87], 
if one could prove MV-ergodicity for these inclusions, then Connes' bicentralizer conjecture would hold true. 
While Theorem 6.4 doesn't bring any progress towards proving MV, MASA, or $R$-ergodicity of such $M \subset \Cal M$, we notice \footnote{I am grateful 
to Stefaan Vaes for his help with the proof of $1^\circ \Rightarrow 4^\circ$ in Corollary 6.7} 
the following equivalences: 

\proclaim{6.7. Corollary} Let $\Cal M$ be a type $\text{\rm III}_1$ factor with $M \subset \Cal M$ the $\text{\rm II}_\infty$ core of its 
continuous decomposition. Let $p\in M$  be non-zero projection  and denote by 
$M_p = pMp \subset p\Cal M p \simeq \Cal M$ the reduced of the inclusion by $p$. 
Consider the conditions: 

\vskip.05in 
$1^\circ$ $\Cal M$ satisfies Connes bicentralizer property. 

\vskip.05in 
$2^\circ$ $M_p \subset \Cal M$ is MV-ergodic.

\vskip.05in 
$3^\circ$ $M_p \subset \Cal M$ is $R$-ergodic

\vskip.05in 
$4^\circ$ $M_p \subset \Cal M$ is MASA-ergodic

\vskip.05in 

Then these conditions are equivalent. Moreover, if one of the conditions $2^\circ-4^\circ$ holds true for some $p\in \Cal P(M)$, then 
they  all hold true for any $p \in \Cal P(M)$.  
\endproclaim
\noindent
{\it Proof}. By Theorem 6.4 we have $4^\circ \Rightarrow 3^\circ$, while $3^\circ \Rightarrow 2^\circ$ 
 is trivial and $2^\circ \Rightarrow 1^\circ$ is (Theorem 3.1 in [H87]).

To see that  $1^\circ \Rightarrow  4^\circ$ for any $p$, note  
that by Haagerup's Theorem ([H87]) if $\Cal M$ satisfies Connes bicentralizer 
property then $\Cal M$ has a normal faithful state $\varphi$ whose centralizer satisfies $\Cal M_\varphi'\cap \Cal M \subset \Cal M_\varphi$.  
By [P81a], this implies $\Cal M_\varphi$ contains a MASA $A_0$ of $\Cal M$. Let $\theta$ be the the 
modular group on $\Cal M$ associated with $\varphi$. Thus, by Takesaki duality ([T74]), the inclusion $M \subset \Cal M$, 
which is isomorphic to its tensor product with $\Cal B(L^2\Bbb R)$, can be viewed as the 
inclusion $\Cal M \rtimes_\theta \Bbb R \subset \Cal M \otimes B(L^2\Bbb R)$.  

In this latter representation of the inclusion $M \subset \Cal M$,  the abelian subalgebra 
$A_0\subset \Cal M \subset M \overline{\otimes} \Cal B(L^2\Bbb R)$ together with $\Bbb R$ generate an abelian von Neumann subalgebra   
$A\subset  M \overline{\otimes} \Cal B(L^2\Bbb R)$ which is generated 
by finite projections of $M\overline{\otimes} \Cal B(L^2\Bbb R)$ and is a MASA in $\Cal M \overline{\otimes} \Cal B(L^2\Bbb R)$, showing that $M \subset \Cal M$ is MASA-ergodic. 
Moreover, when reducing by a projection $p\in A$ that's finite in $M$, 
one gets an abelian algebra in $Ap\subset M_p$ 
which is a MASA in $p\Cal M p\simeq \Cal M$, i.e., $M_p \subset \Cal M$ is MASA-ergodic. 

\hfill $\square$

\heading 7. Some comments and open problems 
\endheading 

It is somewhat frustrating that Theorem 3.3 produces examples of inclusions of factors that are 
MV-ergodic but not $R$-ergodic only ``up to an alternative'' between two 
examples...  Deciding weather a [DPe19]-Poisson boundary inclusion associated with a II$_1$ factor, $M \subset \Cal B_\varphi$   
(which is always MV-ergodic by the Das-Peterson double ergodicity theorem),   
is $R$-ergodic (resp. MASA-ergodic) or not, comes across as a very compelling, important problem that needs clarified. It seems to us that in case $M$ is a 
free group factor, like in Theorem 3.4, then such an inclusion cannot be $R$-ergodic. Note that this would indeed be the case if one could solve in 
the affirmative the following: 

\vskip.05in

\noindent
{\bf 7.1. Conjecture}. {\it The interpolated free group factors $M=L(\Bbb F_t)$, $1<t \leq \infty$, 
contain no ergodic hyperfinite subfactors $L\subset M$ such that $M \subset \langle M, e_L \rangle$  is MV-ergodic}.  
 
 \vskip.05in

It remains as an interesting problem to construct classes of inclusions of separable factors $M \subset \Cal M$ that are MV-ergodic but not 
$R$-ergodic (or if one prefers, not MASA-ergodic).\footnote{Note that the separability of the factors $M, \Cal M$ is essential. Indeed, 
if $M$ is a non-Gamma II$_1$ factor and we let $\Cal M=M^\omega$ for some free ultrafilter $\omega$ then $M'\cap \Cal M=\Bbb C$, implying $M\subset \Cal M$ is MV-ergodic, 
while by [P81a] $M$ contains no MASAs of $\Cal M$, nor 
copies of $R$ that are ergodic in $\Cal M$.}

\vskip.05in

\noindent
{\bf 7.2. Problem}. {\it Identify classes of MV-ergodic inclusions that are not $R$-ergodic} ({\it or MASA-ergodic}).  {\it Alternatively, 
prove that for certain classes of inclusions $($notably $\text{\rm II}_1\subset \text{\rm II}_\infty)$ MV-ergodicity implies $R$-ergodicity, or MASA-ergodicity}.

\vskip.05in

 As we emphasized before, the difficulty in proving the implication in the last part of the Problem 7.2 comes from the bad convexity of the 
$L^1$-norm (as opposed to a Hilbert-norm). 

The class of ergodic II$_1 \subset \text{\rm II}_\infty$ embeddings that one would like to understand most are the ones obtained from the  
basic construction $N \subset M \subset \Cal M=\langle M, e_N \rangle$, where $N\subset M$ is a II$_1$ subfactor with infinite Jones index 
and trivial relative commutant. In this case, $R$-ergodicity amounts to an embedding $R \hookrightarrow M$ such that the Hilbert-bimodule $_RL^2(M)_N$ is irreducible, 
or equivalently $R \vee N^{op}=\Cal B(L^2M)$. This means 
that $\text{\rm sp}RN$ is $\| \ \|_2$-dense in $M$ and more generally $\text{\rm sp}R\xi N$ is dense in $L^2M$ for any non-zero $\xi \in L^2M$. 
One can view this as  $N$ having a  {\it  tight hyperfinite complement} in $M$, an interesting structural property for $N \subset M$. 

Note that by Theorem 3.2, if $M=N \rtimes \Gamma$ for some free action of a group $\Gamma$ on a II$_1$ factor $N$ 
(e.g.,  $N \simeq R$), then 
$N$ has a tight hyperfinite complement in $M$ if and only if $\Gamma$ is amenable. 
Nevertheless, such a crossed product II$_1$ factor $M$ may admit another embedding of $R$ which does have a tight hyperfinite complement.

Before giving such examples, let us fix some terminology.  
If a II$_1$ factor $M$ contains a pair of hyperfinite subfactors $R_0, R_1 \subset M$ that are tight complements one to another, 
then we say that $M$ is $R$-{\it tight}. Thus, $M$ is $R$-tight if it contains an irreducible hyperfinite subfactor $R_1 \subset M$ such that $M \subset \langle M, e_{R_1} \rangle$ is $R$-ergodic. 
Equivalently, there exists a pair of hyperfinite subfactors  $R_0, R_1 \subset M$ such that $R_0 \vee R_1^{op}=\Cal B(L^2M)$. 

A class of examples of tight factors is given by the symmetric enveloping (SE) construction 
in subfactor theory, as introduced in ([P94]). Thus,  if $T \subset S$ denotes 
the SE inclusion of II$_1$ factors arising from an irreducible hyperfinite 
subfactor with finite Jones index, then $S$  is tight (see [P94], [GP96], where the terminology used is ``strongly thin''). Concrete such examples 
are given by crossed product factors $R \rtimes \Gamma$ where $\Gamma \curvearrowright R$ is a free cocycle action of a finitely generated  
group that can split as a diagonal product of two free cocycle actions\footnote{Question: does any free cocycle $\Gamma$-action on $R$ have this 
property?} (for instance, the Bernoulli $\Gamma$-action with base $R$ has this property).

So, despite the fact that $\Gamma$ non-amenable implies the embedding of  $M=N \rtimes \Gamma$ into $\Cal M=\langle M, e_N\rangle$ is not  even MV-ergodic, 
there does exist a hyperfinite subfactor $R_1\subset M$ such that the embedding of $M$ in $\langle M, e_{R_1}\rangle$  is $R$-ergodic. By choosing appropriate 
actions $SL(3, \Bbb Z) \curvearrowright R$ (as in [Cho86]), one can even get 
$M= R \rtimes \Gamma$ to have property (T) (in the sense of [CJ85]). In other words, there are property (T) II$_1$ factors that are $R$-tight. 
  
In view of the iterative technique  of constructing pairs of hyperfinite embeddings into a II$_1$ factor, $R_0, R_1 \subset M$,  
that satisfy given bimodularity properties ([P18]), one would hope that II$_1$ factors satisfying certain  structural properties can be shown to be tight. 
It was speculated in (Conjecture 5.1 in [P18]) that if $M$ is {\it stably single generated} ({\it SSG}) (i.e., there exists $t_n \searrow 0$ such that $M^{t_n}$ is single generated, $\forall n$), 
then $M$ follows tight. As explained in (Section 5 of [P18]), an affirmative answer to this {\it tightness conjecture} 
would imply that $L(\Bbb F_\infty)$ cannot be finitely generated, and the free group factors would follow non-isomorphic (see also [P19] for more on this conjecture). 

An intermediate step towards proving tightness of SSG factors would be to show that such a factor contains an ergodic embedding of a hyperfinite factor 
$R_1\subset M$ such that $M \subset \langle M, e_{R_1} \rangle$ is MV-ergodic (in fact, one would ideally like it to be stably MV-ergodic).

\vskip.05in

\noindent
{\bf 7.3. Problem}. {\it Find sufficient conditions for a} II$_1$ {\it factor $M$ to contain a hyperfinite subfactor $R_1\subset M$ such that $M\subset \langle M, e_{R_1} \rangle$ 
is (stably) MV-ergodic. Do SSG factors satisfy this property} ?

\vskip.05in

Alternatively, in order to show that a II$_1$ factor $M$ is tight, one can try to construct at the same time the two embeddings $R_0, R_1 \subset M$ 
of the hyperfinite II$_1$ factor, using the iterative strategy, so that to have $R_0 \vee R_1^{op}=\Cal B(L^2M)$. A necessary condition 
for the existence of such a ``tight decomposition'' of $M$, indeed an unavoidable  starting point  for constructing $R_0, R_1\subset M$ with $R_0 \vee R_1^{op}=\Cal B(L^2M)$,  
is to have a mean value property of the left-right action of $\Cal U(M)$ on $\Cal B(L^2M)$. 

This motivated us to introduce in an initial version 
of this paper the {\it MV-property} of $M$,  which requires that  
the weak closure of the convex hull of $uv^{op}T{v^{op}}^*u^*$, with $u, v$ running over  $\Cal U(M)$, intersects the scalars, for any $T \in \Cal B(L^2M)$. It also motivated us to ask 
the following questions:

\vskip.05in

\noindent
{\bf 7.4. Problem}.  {\it Do free group factors have the MV-property}? {\it Do all stable single generated factors have the MV-property}? 

\vskip.05in

This problem was meanwhile answered by Das and Peterson in [DPe19], as a consequence of their {\it double ergodicity} theorem. We already 
described this result in Section 3 (see paragraphs preceding Theorem 3.4), but let us recall it here again, in a simplified form 
that better relates to the Problems 7.4 above. Let $M$ be a finitely generated II$_1$ factor, with $\Cal U=\{u_i\}_{i=1}^m\subset \Cal U(M)$ 
a self-adjoint set of  unitary elements with $\Cal U''=M$. For each $T\in \Cal B=\Cal B(L^2M)$ denote  $\varphi(T)=\frac{1}{m}\sum_{i=1}^m u_iTu_i^*$ 
(resp.  $\varphi^{op}(T)= \frac{1}{m}\sum_{i=1}^m u^{op}_iT{u_i^{op}}^*$) the averaging of $T$ by the unitaries in $\Cal U$ (resp. in $\Cal U^{op}=J\Cal UJ$). 
The Das-Peterson double ergodicity theorem states that the averaging of any $T\in \Cal B$ by $\phi=(\varphi + \varphi^{op})/2$ and its powers, 
intersects the scalars, i.e., $\lim_{n\rightarrow \omega} \frac{1}{n} \sum_{k=1}^n \phi^k(T)\in \Bbb C1$, $\forall T\in \Cal B$, where $\omega$ is some fixed free ultrafilter on $\Bbb N$.  
More generally, the same holds true when $\Cal U$ is countable and the averaging in $\varphi$ is taken with positive ``weights'' $\alpha=(\alpha_i)_i$ summing up to $1$. 

So, more than showing that ANY separable (equivalently countably generated) II$_1$ factor has the MV-property, the Das-Peterson theorem provides explicit identification 
of the corresponding left-right averagings. 
This  is quite interesting for  the tightness conjecture, both due to the concreteness of the statement, but also because it shows that if the conjecture is to hold true, 
then the SSG property needs to be used only in the ``second part'' of a potential proof, as described above.

The way we derived from [P81a] that any separable II$_1$ factor $M$ embeds 
ergodically into $R^\infty$ (cf. Corollary 6.6 above) was by first taking an  ergodic embedding of $R$ into $M$ and then using the basic construction 
$M \subset \langle M, e_R \rangle \simeq R^\infty$. It is for such II$_1 \subset \text{\rm II}_\infty$ 
basic construction embeddings that MV-ergodicity, MASA-ergodicity and $R$-ergodicity are most interesting to get. 
But one can ask the similar question in its full abstraction as well: 

\vskip.05in 
\noindent
{\bf 7.5. Problem}. {\it Characterize the class of separable} II$_1$ {\it factors that can be embedded MV-ergodically ({\it respectively MASA/$R$-ergodically}) {\it into} $R^\infty$}.

\vskip.05in

Note that  any ergodic II$_1 \subset \text{\rm II}_\infty$ inclusion, as well as the  II$_\infty$ $\subset \text{\rm III}_1$ inclusions coming from a continuous III$_1$ 
decomposition (together with the $\text{\rm II}_1\subset \text{\rm III}_1$ obtained by reducing them by projections), 
both of which are our main focus of interest, have no normal conditional expectations. 
Nevertheless, inclusions with normal expectations are interesting to study as well. Related to this, it was recently shown 
in [Ma19] that any ergodic inclusion with normal expectation $M \subset \Cal M$ is automatically MV-ergodic. Since  MASA-ergodicity always implies $R$-ergodicity, 
in order to show that all notions of ergodicity, MV, MASA and $R$-ergodicity, coincide for such inclusions, one only needs 
to answer:

\vskip.05in 

\noindent
{\bf 7.6. Problem} {\it Does  MV-ergodicity imply MASA-ergodicity  for inclusions of factors with normal expectation ?}

\head  References \endhead

\item{[AP17]} C. Anantharaman, S. Popa: ``An introduction to II$_1$ factors'', \newline www.math.ucla.edu/$\sim$popa/Books/IIun-v13.pdf

\item{[ChHo08]} I. Chifan, C. Houdayer: {\it Bass-Serre rigidity results in von Neumann algebras}, Duke Math. J. {\bf 153} (2010), 23-54. 

\item{[Cho86]} M. Choda: {\it Outer actions of groups with property} (T) {\it on the hyperfinite} II$_1$ {\it factor}, Math. Japan. {\bf 31}, 
(1986), 533-551.

\item{[C73]} A. Connes: {\it Une classification des facteurs de type} III, Annales de l'Ec. Norm Sup., {\bf 6} (1973), 133-252. 

\item{[C75]} A. Connes: {\it On the classification of von Neumann algebras and their automorphisms}, 
Symposia Mathematica, Vol. XX (Convegno sulle Algebre C$^*$, INDAM, Rome, 1975), Academic Press, London, 1976, p. 435-478.

\item{[C76]} A. Connes: {\it Classification of injective factors}, Ann. of Math., {\bf 104} (1976), 73-115. 

\item{[C85]} A. Connes: {\it Type} III$_1$ {\it factors, property} $L_\lambda'$ 
{\it and closure of inner automorphisms},  J. Operator Theory, {\bf 14} (1985), 189-211.

\item{[CJ85]} A. Connes, V.F.R. Jones: {\it Property} (T) {\it for von Neumann algebras}, Bull. London Nath. Soc. {\bf 17} (1985), 57-62. 

\item{[CT77]} A. Connes, M. Takesaki: {\it The flow of weights of factors of type} III, Tohoku Math. Journ. {\bf 29} (1977), 473-575.

\item{[DPe19]} S. Das, J. Peterson: {\it Poisson boundaries of} II$_1$ {\it factors}, preprint 2019. 

\item{[D57]} J. Dixmier: ``Les alg\'ebres d'op\'erateurs dans l'espace Hilbertien,'' Gauthier-Villars, Paris, 1957.

\item{[Dy93]} K. Dykema: {\it Interpolated free group factors}, Pacific J. Math. 163 (1994), 123-135.

\item{[FK51]} B. Fuglede,  R. V. Kadison: {\it On a conjecture of Murray and von Neumann}, Proc. Nat.
Acad. Sci. U.S.A. {\bf 36} (1951), 420-425.

\item{[GP98]} L. Ge, S. Popa: {\it On some decomposition properties for factors of type} II$_1$, Duke
Math. J. {\bf 94} (1998), 79-101.

\item{[H87]} U. Haagerup: {\it Connes' bicentralizer problem and uniqueness of the injective factor of type} III$_1$, 
Acta. Math., {\bf 158} (1987), 95-148. 

\item{[HoP16]} C. Houdayer, S. Popa: {\it Singular MASAs in type} III {\it factors and Connes' bicentralizer property}, 
Proceedings of the 9th MSJ-SI "Operator Algebras and Mathematical Physics" held in Sendai, Japan (2016), math.OA/1704.07255.   

\item{[J83]} V.F.R. Jones: {\it Index for subfactors}, Invent. Math. {\bf 72} (1983), 1-25. 

\item{[K67]} R.V. Kadison: {\it Problems on von Neumann algebras.} Baton Rouge Conference, 1967 (unpublished).

\item{[L84]} R. Longo: {\it Solution of the factorial Stone-Weierstrass conjecture}, Invent. Math. {bf 76} (1984), 145-155. 

\item{[Ma19]} A. Marrakchi: {\it On the weak relative Dixmier property}, arXiv:1907.00615

\item{[MvN36]} F. Murray, J. von Neumann: {\it On rings of operators}, Ann. Math. {\bf 37} ,  (1936),  116-229. 

\item{[MvN43]} F. Murray, J. von Neumann: {\it On rings of operators  IV}, Ann. Math. {\bf 44} (1943), 716-808.

\item{[vN29]} J. von Neumann: {\it Zur Algebra der Funktionaloperatoren und Theorie der 
normalen Operatoren}, Math. Ann. {\bf 102} (1929), 370-427

\item{[vN54]} J. von Neumann: {\it Unsolved problems in mathematics}, Plenary Talk,  ICM 1954, 
 in  ``John von Neumann and the foundation of quantum physics'', Vienna Circle Institute Yearbook [2000], {\bf 8}, pp 226-248. 

\item{[P81a]} S. Popa: {\it On a problem of R.V. Kadison on maximal abelian *-subalgebras in factors}, Invent. Math., {\bf 65} (1981),
269-281.

\item{[P81b]} S. Popa: {\it Maximal injective subalgebras in factors associated with free groups}, Advances in Math., {\bf 50} (1983), 27-48. 

\item{[P81c]} S. Popa: {\it Singular maximal abelian *-subalgebras in
continuous von Neumann algebras}, J. Funct. Analysis, {\bf 50} (1983), 151-166,

\item{[P83]} S. Popa: {\it Hyperfinite subalgebras normalized by a
given automorphism and related problems}, in ``Proceedings of the
Conference in Op. Alg. and Erg. Theory'' Busteni 1983, Lect. Notes
in Math., Springer-Verlag, {\bf 1132}, 1984, pp 421-433.

\item{[P84]} S. Popa, {\it Semiregular maximal abelian *-subalgebras
and the solution to the factor state Stone-Weierstrass problem},
Invent. Math., {\bf 76} (1984), 157-161.

\item{[P94]} S. Popa: {\it Symmetric enveloping algebras,
amenability and AFD properties for subfactors}, Math. Res.
Letters, {\bf 1} (1994), 409-425.

\item{[P98]} S. Popa: {\it On the relative Dixmier property for inclusions of} C$^*$-{\it algebras}, Journal of Funct. Analysis, {\bf 171} 
(2000), 139-154.

\item{[P01]} S. Popa: {\it On a class of type} II$_1$ {\it factors with
Betti numbers invariants}, Ann. of Math {\bf 163} (2006), 809-899
(math.OA/0209310; MSRI preprint 2001-024).

\item{[P17]} S. Popa:  {\it  Asymptotic orthogonalization of subalgebras in} II$_1$ {\it factors}, Publ. RIMS Kyoto University, {\bf 55} (2019), 795-809 (math.OA/1707.07317)

\item{[P18a]} S. Popa: {\it On the vanishing cohomology problem for cocycle actions of groups on}  II$_1$ {\it factors}, 
to appear in Ann. Ec. Norm. Super., math.OA/1802.09964

\item{[P18b]} S. Popa: {\it Coarse decomposition of} II$_1$ {\it factors}, math.OA/1811.09213 

\item{[P19]} S. Popa: {\it Tight decomposition of factors and the single generation problem}, to appear in J. Operator Theory 
vol. dedicated to D. Voiculescu's 70th anniversary. 

\item{[PV14]} S. Popa, S. Vaes: {\it Vanishing of the continuous first L$^2$-cohomology for} II$_1$ {\it factors}, 
Intern. Math. Res. Notices {\bf 8} (2014), 1-9 (math.OA/1401.1186)

\item{[Po67]} R. Powers: {\it Representations of uniformly hyperfinite algebras and their associated von Neumann rings}, Ann. Math. {\bf 86} (1967), 
138-171. 

\item{[R92]} F. Radulescu: {\it Random matrices, amalgamated free products and subfactors of the von Neumann algebra of a free group, of noninteger 
index}, Invent. Math. {\bf 115} (1994), 347-389.

\item{[S71]} S. Sakai: ``C$^*$-algebras and W$^*$-algebras'', 1971, Springer Verlag.  

\item{[Sc63]} J. Schwartz: {\it Two finite, non-hyperfinite, non-isomorphic factors}, Comm. Pure App. Math. (1963), 19-26. 

\item{[Sk77]} C. Skau: {\it Finite subalgebras of a von Neumann algebra}, J. Funct. Anal. {\bf 25} (1977), 211-235. 

\item{[T74]} M. Takesaki: {\it Duality for crossed products and the structure of yon Neumann algebras of type} III, 
Acta Math., {\bf 131} (1974), 249-310.

\item{[T03]} M. Takesaki: ``Theory of operator algebras II '', 
Encyclopaedia of Mathematical Sciences, {\bf 125}. 
Operator Algebras and Non-commutative Geometry, 6. Springer-Verlag, Berlin, 2003. xxii+518 pp.

\enddocument